\documentclass{article}
\usepackage{latexsym,amssymb}
\usepackage[dvips]{graphicx}
\usepackage{subfigure}

\newcommand{\C}{\mathbb{C}}
\newcommand{\R}{\mathbb{R}}
\newcommand{\T}{\mathbb{T}}
\newcommand{\I}{\mathbb{I}}

\newcommand{\oct}{\mathbb{O}}
\newcommand{\<}{\left<}
\newcommand{\m}{\right>}

\newcommand{\oo}{\mathcal{O}}
\newcommand{\p}{\mathcal{P}}
\newcommand{\loup}{\mathcal{L}}
\newcommand{\J}{J}
\newcommand{\Fix}{\mathop{\rm Fix}}
\newcommand{\Ker}{\mathop{\rm Ker}}
\newcommand{\x}{x}

\newcommand{\so}{\mathfrak{so}(3)}
\newcommand{\se}{\mathfrak{se}(2)}
\newcommand{\sodeux}{\mathfrak{so}(2)}
\newcommand{\g}{\mathfrak{g}}
\newcommand{\liel}{\mathfrak{l}}

\def\Z{{\sf Z\hspace{-.45em}Z}}

\def\rmq{\noindent\textbf{Remark. }}

\def\1N{1\dots N}

\def\eq{relative equilibrium }
\def\eqb{relative equilibrium}
\def\eqs{relative equilibria }
\def\eqsb{relative equilibria}
\def\ops{periodic orbits }
\def\opsb{periodic orbits}

\def\rum{relative equilibrium }
\def\rumb{relative equilibrium}
\def\ria{relative equilibria }
\def\riab{relative equilibria}
\def\rmk{\noindent\textbf{Remark. }}

\def\mom{momentum map }
\newtheorem{theo}{Theorem}[section]

\newtheorem{prop}[theo]{Proposition}

\newenvironment{demo}{\addvspace\baselineskip\noindent{\textbf{Proof.}}\quad}{\hspace*{\fill} $\Box$\par\addvspace\baselineskip}
\newenvironment{proof}{\addvspace\baselineskip\noindent{\textbf{Proof.}}\quad}{\hspace*{\fill} $\Box$\par\addvspace\baselineskip}

\def\hfl#1#2{\smash{\mathop{\hbox to 12mm{\rightarrowfill}}\limits^{\scriptstyle#1}_{\scriptstyle#2}}}

\def\CC{C}
\def\DD{D}

\def\kO{{\cal O}}
\def\OO{O}

\def\SO{SO}
\def\SS{{\cal S}}

\def\restr#1{\vrule height1.2ex width.4pt
               depth1.4ex\lower1.0ex\hbox{\scriptsize $\,#1$}}

\def\Tet{\mathbb{T}}
\def\Oct{\mathbb{O}}
\def\Icos{\mathbb{I}}

\newtheorem{remark1}[theo]{Remark}
\newtheorem{example1}[theo]{Example}

\def\ie{that is }
\newcommand{\bfig}{\begin{figure}[htbp]}
\newcommand{\efig}{\end{figure}}
\newcommand{\bc}{\begin{center}}
\newcommand{\ec}{\end{center}}
\newcommand{\beq}{\begin{equation}}
\newcommand{\eeq}{\end{equation}}
\newcommand{\pg}{\textquotedblleft}
\newcommand{\pd}{\textquotedblright\ }

\begin{document}

\noindent{\Large\textbf{Relative periodic orbits in point vortex systems}}\\
\\
F Laurent-Polz\\
Institut Non Lin\'eaire de Nice, Universit\'e de Nice - Sophia Antipolis, Valbonne, France\\

\noindent\textbf{Abstract.}
We give a method to determine relative periodic orbits in point vortex systems:
it consists mainly into perform a symplectic reduction on a fixed point submanifold in order
to obtain a two-dimensional reduced phase space.
The method is applied to point vortices systems on a sphere and on the plane, but works
for other surfaces with isotropy (cylinder, ellipsoid, \dots).
The method permits also to determine some relative equilibria and heteroclinic cycles
connecting these relative equilibria.\\


\noindent AMS classification scheme number: 58Z05, 70H14, 70H33\\

\section{Introduction}
\label{intro}

Since the work of Helmholtz \cite{H}, systems of point vortices on a surface
have been studied as finite-dimensional approximations to vorticity dynamics in an ideal fluid.
They model the dynamics of concentrated regions of vorticity, such as cyclones and hurricanes.
A parameter $\lambda$ is attached to a point vortex: its \emph{vorticity},
it is the value of the vorticity at the locus of the vortex,
or equivalently the value of the fluid circulation around the vortex.
For a general survey on point vortices, see \cite{ANSTV03}.
In this paper we will consider point vortex systems on both the sphere and the plane,
which are Hamiltonian systems with symmetry.
Concerning the spherical case, we consider a \emph{non-rotating} sphere since the rotation of the sphere
induces a non-uniform background vorticity which makes the whole system infinite-dimensional.

Relative equilibria are orbits of the symmetry group action which are invariant under the flow, this
corresponds here to motions of point vortices which are stationary in a steadily rotating frame.
We first review the literature about the point vortex system on the sphere.
The existence and nonlinear stability of relative equilibria formed of three vortices have been studied
respectively in \cite{KN98} and \cite{PM98}.
Existence of relative equilibria formed of $N$
vortices is treated in \cite{LMR00}, and the nonlinear stability of a
latitudinal ring of $N$ identical vorticities is computed in \cite{BC01,LMR}.
Existence and nonlinear stability of relative equilibria of $N$ vortices of vorticity $+1$ together with
$N$ vortices of vorticity $-1$ are studied in \cite{LP02}.
The existence and nonlinear stability of more complex arrangements, such as two latitudinal rings of
identical vortices with and without polar vortices, are studied in \cite{LMR}.
It has also been proved in \cite{LP} that any relative equilibrium formed of
latitudinal rings of identical vortices for the non-rotating sphere
persists to be a relative equilibrium when the sphere rotates.

Relative periodic orbits (RPOs for short) are the analogous of relative equilibria concerning periodic orbits,
this corresponds here to motions which are periodic in a steadily rotating frame
(a precise definition is given in Section \ref{methode}).
Periodic orbits on the sphere were determined in \cite{ST,To} thanks to the following method:
they reduced the system to two-dimensional systems by a symmetric reduction
(using some finite subgroups of $SO(3)$); the computation of the dynamics on the reduced spaces
permits then to determine periodic orbits.
Our paper is devoted to transpose that method to determine \emph{relative} periodic orbits.
To this end, we combine a symmetric reduction together with a symplectic reduction.
The method is explained in details in Section \ref{methode}, it permits to determine
periodic orbits, RPOs, and heteroclinic cycles connecting relative equilibria,
this for an \emph{arbitrary} number of vortex.

Section \ref{rposphere} is devoted to the application of the previous method to point vortices systems on a
sphere. We recall the symmetries and conserved quantities of that system, and show that four
finite groups of $O(3)$ ($C_n,C_{nh},C_h$ and $C_i$) permit to apply the method.
The groups $C_n$ and $C_{nh}$ give RPOs (as well as relative equilibria, periodic orbits,
heteroclinic cycles between relative equilibria) for an arbitrary number of vortices.
The group $C_i$ gives RPOs formed of up eight vortices.
We apply the method first for systems of $k_r N+k_p$ vortices ($k_r,k_p=0,1$ or $2$) where
the $N$ first (resp. the $N$ second) vortices are identical (that is have equal vorticities).
We obtain RPOs formed of two rings of each $N$ identical vortices (a $N$-ring) with and without one or two
polar vortices.
Among the previous systems, there is a particular one which merits further attention.
Indeed if we consider systems of $N$ vortices of vorticity $+1$ and
$N$ vortices of vorticity $-1$, we obtain more symplectic symmetries.
This permits to determine RPOs formed of up to four $N$-rings,
and interesting heteroclinic cycles between relative equilibria.
Some bifurcations are also described.

Section \ref{rpoplane} is devoted to the application of the method to point vortices systems on the plane.
Numerous papers have been written on vortices in the plane \cite{Ar82,Ar83',Ar98,LR96}.
The nonlinear stability of the $C_N(R)$ relative equilibria (a planar $N$-ring) has been
calculated in \cite{CS99} a hundred years after the early works of Thomson.
The nonlinear stability of the $C_N(R,p)$ relative equilibria (a planar $N$-ring plus a central vortex)
is calculated also in \cite{CS99} (see also the Appendix of \cite{LP}).
In order to apply the method, the planar problem is less rich than the spherical one
since we dispose of only two types of finite groups: $C_n$ and $D_n$.
Due to the noncompactness of the translational symmetries, the method could not determine
relative equilibria and RPOs with a general translational motion.
We obtain RPOs with a general rotational motion, and we success to obtain some unbounded
motions. A heteroclinic cycles between relative equilibria is also determined.

\section{Method for determining relative periodic orbits}
\label{methode}

We first review some generalities about Hamiltonian systems with symmetry, and relative periodic
orbits.

Let $(\p,\omega,H,G)$ a Hamiltonian system with symmetry such that the action of $G$ on $\p$ is
semi-symplectic, that is:

$\bullet$ $(\p,\omega)$ is a symplectic manifold,

$\bullet$ $G$ is a Lie group acting smoothly on $\p$,

$\bullet$ $H: \p\to\R$ is $G$-invariant function (the Hamiltonian),

$\bullet$ the action of $G$ is semi-symplectic ($g^*\omega=\pm\omega$ for all $g\in G$).

\noindent The Hamiltonian vector field $X_H$ is defined by $\omega(X_H,\cdot)=dH$,
and defines a dynamical system $\dot v=X_H(v)$.
An element $g\in G$ is said to be \textit{symplectic} (resp. \textit{anti-symplectic})
if it preserves the symplectic form (resp. changes into its opposite),
that is if $g^*\omega=\omega$ (resp. $g^*\omega=-\omega$).
The vector field $X_H$ is $G^o$-equivariant where $G^o$ is
the connected component of identity in $G$ (which is a group formed of symplectic elements).
The fixed point set of a subgroup $K$ of $G$ is
$\Fix (K,\p)=\lbrace x\in\p \mid g\cdot x = x, \forall g\in K  \rbrace$
which is a closed submanifold of $\p$;
the normalizer of $K$, $N_G(K)=\lbrace g\in G \mid g\cdot K \cdot g^{-1} = K \rbrace$,
acts on $\Fix (K,\p)$ by restriction.
If $K<G$ is formed of symplectic elements, then $X_H$ is $K$-equivariant and $\Fix (K,\p)$
is invariant by the flow of the dynamical system.
If in addition the action of $K$ on $\p$ is proper, then $\Fix (K,\p)$ is a Hamiltonian subsystem
with Hamiltonian given by the restriction of $H$ to $\Fix (K,\p)$: this is called \emph{symmetric reduction}.
We will write $\Fix K$ instead of $\Fix (K,\p)$ when there are no ambiguities.

Let $F_t$ the flow of $X_H$. A point $p\in\p$ is said to be \emph{periodic} if there exists a
constant $T>0$ such that for all time $t$, $F_{t+T}(p)=F_t(p)$. The \emph{period} is the smallest
$T>0$ which satisfies that condition. The set $\gamma=\{F_t(p)\mid t\in\R\}$ is called a
\emph{periodic orbit}. Every point of $\gamma$ is periodic with the same period, hence we can
define \emph{the} period of a periodic orbit.

Relative equilibria are dynamical trajectories that are generated by
the action of a continuous connected subgroup of the symmetry group:
for all $t$ there exists $g_t \in G^o$ such that $x(t)=g_t \cdot x(0)$.
More intuitively, this will correspond here to motions of the point vortices
which are stationary in a steadily rotating frame; in particular these motions
are here periodic orbits.
In other words, the motion of a \rum corresponds here to a rigid rotation of $N$ point
vortices about some axis (resp. around some point) in $\R^3$ (resp. in in $\R^2$).

A point $p\in\p$ is said to be a \emph{relative periodic point} if there exist $g\in G^o$ and
$T>0$ such that $F_{t+T}(p)=g\cdot F_t(p)$ for all time $t$. The set $\gamma=\{F_t(p)\mid
t\in\R\}$ is called a \emph{relative periodic orbit} (RPO), and every point of $\gamma$ is a relative
periodic point. In particular, a periodic orbit which is not a \eq is a RPO. Note that there exist
other equivalent definitions of a relative periodic orbit \cite{Or98}.
Typically a RPO is a solution which, in a suitably moving frame, looks time-periodic.

Recall two well-known results about periodic orbits which have been recently generalized
to RPOs. In Hamiltonian systems with symmetry, RPOs are
typically present around stable \eqs \cite{LT99,OR}, and, they persist
when changing the values of both the Hamiltonian and the momentum
(provided the RPO is non-degenerated and the action of $G^o$ is free) \cite{Mo97b}.
Splitting the dynamics around a RPO permits to yield a local theory about stability and persistence
of RPOs in Hamiltonian systems with symmetry \cite{WR}.
We will not go further since we are here only interested by existence of RPOs.


\bigskip

The idea in order to determine \ops is to find isotropy subgroups $K\subset G$ such
that:
\begin{itemize}
  \item $K$ is formed of symplectic elements,
  \item the action of $K$ on $\p$ is proper (this holds if $K$ is compact),
  \item $\Fix K$ is two-dimensional.
\end{itemize}
Indeed if these three conditions hold, then $\Fix K$ is a Hamiltonian  sub-system with
Hamiltonian given by the restriction of $H$ to $\Fix K$. Moreover, there
exist \ops in a two-dimensional Hamiltonian system for which the
Hamiltonian is a non-constant proper map, as the following proposition shows.
These \ops can be easily drawn in the phase space by calculating the level sets of $H|_{\Fix K}$,
since $\Fix K$ is two-dimensional.

One can show that a continuous function $h:M\subset\R^{2n}\to\R$ is proper if $h$
takes infinite values on $\overline{M}\setminus M$ and if $\overline{M}$ is compact.
This means for point vortices that the restriction of $H$ to $\Fix(K,\p)$ is a proper function if the
Hamiltonian takes infinite values on
$\overline{\Fix(K,\p)}\setminus\Fix(K,\p)\subset\Fix(K,\overline{\p})\setminus \Fix(K,\p)$,
that is if there are not possible collision configurations in $\Fix(K,\overline{\p})$
from a dynamical point of view. Indeed if a collision is dynamically possible, then the Hamiltonian
takes a finite value at the collision point.

\begin{prop}
\label{proprenoncste} Let $(M,\omega,h)$ a Hamiltonian system such that the phase space $M$
is two-dimensional and the Hamiltonian $h$ is a proper map. Let $E$ a regular value
of the Hamiltonian $h$, \ie $dh(x)\neq 0$ for all $x\in h^{-1}(E)$
($h^{-1}(E)$ does not contain equilibrium points).
Then the connected components of $h^{-1}(E)$ are trajectories of periodic orbits in $M$.
\end{prop}
\begin{demo}
Let $C_E$ a connected component of $h^{-1}(E)$. Since $E$ is a regular value of $h$ and $\dim
M=2$, the set $h^{-1}(E)$ is a one-dimensional manifold. Moreover, the manifold $h^{-1}(E)$ is
compact since $h$ is proper.
One-dimensional compact manifolds are diffeomorphic to a union of circles,
hence $C_E$ is diffeomorphic to a circle. Moreover, the velocity never vanishes on
$C_E$ since $E$ is a regular value of $h$ ($dh(x)\neq0$ for all $x\in h^{-1}(E)$), thus $C_E$
is the trajectory of a periodic orbit.
\end{demo}

\bigskip

In order to determine RPOs, the method is similar. We look for compact subgroups $K$ formed of
symplectic elements, then $(M\!:=\!\Fix K,\tilde\omega,\tilde H)$ is a Hamiltonian sub-sytem
(the \emph{tilde} denotes the restriction to $M$). Assume that $K\neq\{Id\}$ and there exists
a $G^o$-equivariant momentum map $\J:\p\to\g^*$, assume in addition that the action of $G^o$ on
$\p$ is proper and free, hence we can apply the \emph{symplectic reduction} \cite{MR94}:
the reduced spaces $(\p_\mu=\J^{-1}(\mu)/G_\mu^o,H_\mu,\omega_\mu)$ are Hamiltonian where
$T\varpi^*\omega_\mu=\omega$, $H_\mu\varpi=H$, and $\varpi$ is the projection $\varpi:\p\to\p_\mu$.
We will consider finite symmetry groups $K$, thus the assumption \pg $K$ compact\pd will be
automatically satisfied.
These assumptions hold for point vortices on sphere if $N>2$, and on plane taking $G=O(2)$ only
(not the full symmetry group $E(2)$ with its translational symmetries).

The normalizer of $K$, $N_G(K)$, acts on $M$, but this action is not free
since $K\subset N_G(K)$ and $K\subset G_x$ for all $x\in M$. However, the group $N_G(K)^o$
(the connected component of the identity in $N_G(K)$) acts properly and freely on $M$ since
$N_G(K)^o\subset G^o$. The system $(M,\tilde\omega,\tilde H,N_G(K)^o)$ is a Hamiltonian system
with symmetry, and the properness and freeness of the action allow us to perform a symplectic
reduction. Assume that $\dim N_G(K)^o\neq 0$ and denote $L=N_G(K)^o$ for convenience.
There exist a $L$-equivariant \mom $\J_L:M\to\liel^*$ by the next proposition,
and the reduced systems
$(M_\nu=\J_L^{-1}(\nu)/L_\nu,\tilde\omega_\nu,\tilde H_\nu)$ with $\nu\in\liel^*$ are Hamiltonian.
If the reduced space $M_\nu$ is two-dimensional, then trajectories are generically \ops in $M_\nu$
(Prop. \ref{proprenoncste})
\ie RPOs in $M$ (hence RPOs in $\p$).
We will therefore look for groups $K$ such that $\dim M_\nu=2$.

\begin{prop}
\label{propmom}
Let $K$ be a subgroup of $G$, and $L=N_G(K)^o$.
Assume that there exists a $G$-equivariant \mom $\J:\p\to\g^*$, and that $K$ is formed of
symplectic elements.
Then the map $J_L=p\circ\J\circ i : \Fix(K,\p)\to\liel^*$ is a $L$-equivariant momentum map,
where $i:\Fix(K,\p)\hookrightarrow\p$ and $p:\g^*\to\liel^*$ are respectively the canonical injection
and projection.
\end{prop}

\begin{proof}
The existence comes from the defining equation for the original momentum map
$$
\< d\J_x(v),\xi \m=\omega_x(\xi_\p(x),v),\ \forall x\in\p,\ \forall v\in T_x\p,\ \forall \xi\in\g
$$
restricted to the fixed point (symplectic) manifold $\Fix(K,\p)$.
The $L$-equivariance is due to the $L$-equivariance of $p$ (coadgoint action) together with
the $G$-equivariance of the other two maps $J$ and $i$.
\end{proof}

\vspace{1cm}

\noindent\textbf{{Summary of the method}}

\bigskip\noindent
In order to determine RPOs (with isotropy), the first step consists to find subgroups $K$ of $G$ such that:

\bc
\begin{tabular}{rl}
   $\bullet$ & $K$ is compact,\\
   $\bullet$ & $K\subset\Ker\chi$,\\
   $\bullet$ & $\dim N_G(K)^o\neq 0$,\\
   $\bullet$ & $\dim M_\nu =2$.
\end{tabular}
\ec


\noindent Thanks to the level sets of $\tilde H_\nu$, the phase diagram in $M_\nu$ can then be easily computed,
and we have the following correspondance:


\bc
\begin{tabular}{c|c|c}
   $M_\nu$ & & $\p$\\
   \hline
   \hline
   & & \\
  equilibrium & $\Longrightarrow$ & relative equilibrium\\
  periodic orbit & $\Longrightarrow$ & RPO \\
  heteroclinic orbit & $\Longrightarrow$ & heteroclinic orbit between\\
   & & relative equilibria\\
\end{tabular}
\ec

\bigskip

\rmq This method is useful overall for $N$-body systems:
indeed, in that systems there is a continuous symmetry group $G^o$ for the vector field acting freely,
which permits to perform the symplectic reduction;
but the symmetry group $G$ of the Hamiltonian contains permutations symmetries and does not act
freely, hence we can restrict the dynamics to the fixed points manifolds $\Fix G_x$ (symmetric reduction).

\bigskip

In the following sections, we will study dynamics in the spaces $M$ and $M_\mu$,
the notion of stability  should be understood as stability \emph{in} that spaces.

\section{Point vortices on a sphere}
\label{rposphere}

\subsection{Description of the Hamiltonian system}
\label{pointsphere}

We describe the $N$-vortex system on a sphere which is a $N$-body Hamiltonian system with symmetry,
proofs and details can be found in \cite{KN98,LMR00,LP02}.
Consider $N$ vortices $x_1,\dots,x_{N}$ on the sphere $S^2$ with vorticities $\lambda_1,\dots,\lambda_{N}\in\R$.\\

Let $\theta_i,\phi_i$ be respectively the co-latitude and the longitude of the vortex $x_i$.
The dynamical system is Hamiltonian with Hamiltonian given by
$$
H = \sum_{i<j}\lambda_i\lambda_j\ln ( 1-\cos\theta_i\cos\theta_j-\sin\theta_i\sin\theta_j\cos(\phi_i-\phi_j) )
$$
and conjugate variables given by
$q_i=\sqrt{\vert\lambda_i\vert}\cos\theta_i$ and
$p_i=sign(\lambda_i)\sqrt{\vert\lambda_i\vert}\phi_i$.

The phase space is
$\p = \lbrace (x_1,\dots,x_{n})\in S^2\times\cdots\times S^2\mid  x_i\neq x_j\ \mbox{if}\ i \neq j
\rbrace$ endowed with the symplectic form $\omega=\sum_{i} \lambda_i \sin\theta_i\ d\theta_i\land
d\phi_i$.
The Hamiltonian vector field $X_H$ satisfies $\omega(X_H(x),\cdot)=dH_x$.
If the vortices $x_j\in S^2$ are embedded in $\R^3$, then we obtain
$$
\dot {x_i} = X_H(x)_i=\sum_{j,j\neq i} \frac{\lambda_{j}(x_{j}\times x_{i})}{1-x_{i}\cdot x_{j}}\  ,\
i=1,\dots,N,
$$
$$
H = \sum_{i<j}\lambda_i\lambda_j\ln(\|x_i-x_j\|^2/2).
$$
\noindent\textbf{{Symmetries.}}
It follows from expressions of $H$ and $X_H$ that
$X_H$ is $SO(3)$-equivariant, and $H$ is $O(3)$-invariant.
Moreover $X_H$ and $H$ are respectively equivariant and invariant
under permutations of vortices with equal vorticities.
Let precise that statement.
The group $O(3)\times S_N$ acts on $\p_N$ in the following manner:
$$(A,\sigma)(\x_1,\dots,\x_{N})=(A\x_{\sigma(1)},\dots,A\x_{\sigma(N)}).$$
Let $\loup=(\lambda_1,\dots,\lambda_N)$ and
$$\hat S(\loup)=\{\sigma\in S_N \mid\
\exists\;\varepsilon(\sigma)=\pm 1 ,\ \forall i,\
 \lambda_{\sigma(i)}=\varepsilon(\sigma)\lambda_i \}.$$
$$
\begin{array}{rcl}
\mbox{Let}\ \chi: O(3)\times\hat S(\loup) & \to & \lbrace -1,+1 \rbrace \\
(A,\sigma) & \mapsto & \varepsilon(\sigma)\ \det A
\end{array}
$$
It is then straightforward to verify that:
\emph{The Hamiltonian $H$ is invariant under $G=O(3)\times\hat S(\loup)$,
and the vector field $X_H$ is $\Ker(\chi)$-equivariant.}

For example in the case of $N$ identical vortices (vortices of equal vorticities),
we have $G=O(3)\times S_N$ and $\Ker(\chi)=SO(3)\times S_N$.
In the case of $N$ vortices of vorticity $+1$ and $N$ vortices of vorticity $-1$,
we have $G=O(3)\times S_N\times S_N\rtimes \Z_2 \lbrack \tau \rbrack$
and $\Ker(\chi)=SO(3)\times S_N\times S_N\rtimes \Z_2 \lbrack (-Id,\tau) \rbrack$,
where $\tau=+/- \in S_{2N}$ is the permutation which exchange vortices $(+1)$ with vortices $(-1)$.

\bigskip

\noindent\textbf{{Conserved quantities.}}
The \mom $J : \p_N\to \g^*\!=\!\so^*\!\simeq\!\R^3$ for the diagonal action of $G^o=SO(3)$ on $\p_N$ is:
$$
J(x)=\sum_{j=1}^{N}\lambda_j\; x_j
$$
where $x_j\in S^2\subset\R^3$.
The three components of the \mom are conserved under the dynamics.
We can choose then a frame $(O,\vec{e_x},\vec{e_y},\vec{e_z})$ of $\R^3$
such that $\vec{e_z}=\J/\|\J\|$ (provided $\J\neq0$).

\bigskip

\noindent\textbf{{Application of the method.}}
From Table \ref{normaliz}, the groups $C_n,C_{nh},C_h$ and $C_i$ are the only groups for which
the normalizer has a non-zero dimension.
Thus only these groups will permit to perform the method described in Section \ref{methode}.
\begin{table}[tb]
\begin{center}
\begin{tabular}{|c|c|c|}
\hline $K$ & $N_{SO(3)}(K)$ & $N_{SO(3)}(K)^o$\\ \hline $C_{n}$ & $D_\infty$ & $SO(2)$ \\ $D_2$ &
$\oct$ & 1 \\ $D_n,\ n>2$ & $D_{2n}$ & 1 \\ $\T$ & $\oct$ & 1 \\ $\oct$ & $\oct$ & 1 \\ $\I$ &
$\I$ & 1 \\ \hline \hline $K$ & $N_{O(3)}(K)$ & $N_{O(3)}(K)^o$\\ \hline $C_{i}$ & $SO(3)\times
C_i$ & $SO(3)$ \\ $C_h$ & $D_{\infty h}$ & $SO(2)$ \\ $C_{nh}$ & $D_{\infty h}$ & $SO(2)$ \\
$C_{nv}$ & $D_{2nh}$ & 1 \\ $D_{nh}$ & $D_{2nh}$ & 1 \\ $D_{nd}$ & $D_{2nh}$ & 1 \\ $\T_h$ &
$\oct_h$ & 1 \\ $\T_d$ & $\oct_h$ & 1 \\ $\oct_h$ & $\oct_h$ & 1 \\ $\I_h$ & $\I_h$ & 1 \\ \hline
\end{tabular}
\caption{Normalizers of finite subgroups of $O(3)$.}
\label{normaliz}
\end{center}
\end{table}
For $K=C_n,C_{nh},C_h$, we have $L=N_G(K)^o=SO(2)$.
Moreover, $SO(2)_\nu=SO(2)$ for all $\nu\in\sodeux^*\simeq\R$,
the manifold $M_\nu$ will be of dimension
$\dim\Fix K - 2$.

We then look for cases where $\dim\Fix K=4$ in order to obtain $\dim M_\nu=2$.
We have just three types of groups, we must therefore play with $N$ the number of vortices to get
$\dim\Fix K=4$. Indeed $\Fix K=\Fix(K,\p_N)$, the fixed points manifold depends on $\p_N$.
In our case of point vortices on a sphere,
$\dim\Fix K=4$ will mean that $\Fix K=S^2\times S^2\setminus\Delta$ where
$\Delta$ is the set of possible collisions in $\Fix K$, and this manifold will be parametrized by
coordinates $(\theta_1,\phi_1),(\theta_1^\prime,\phi_1^\prime)$ of two vortex of, say
vorticities $\lambda_1$ and $\lambda_1^\prime$. From Proposition \ref{propmom}, the \mom of
the action of $L=SO(2)$ on $\Fix K$ is given by
$$\J_L= n\cdot(\lambda_1\cos\theta_1+\lambda_1^\prime\cos\theta_1^\prime)$$
and $\J_L$ is $SO(2)$-equivariant.
The restriction $\tilde\J$ of $\J$ to $\Fix K$ satisfies
$\tilde\J(x)=(0,0,\J_L(x))\in\R^3$, hence we identify $\tilde\J$ with its $z$-component $\J_L$.
So for  $\nu=\J_L(x)\in\R$ we will write
$M_\mu$ (resp. $\tilde H_\mu$) instead of $M_\nu$ (resp. $\tilde H_\nu$) where $\mu=\tilde\J(x)\in\R^3$.

For $K=C_i$, we have $L=N_G(K)^o=SO(3)$, and for $\mu\in\so^*\simeq\R^3$, $SO(3)_\mu=SO(2)$
if $\mu\neq 0$ and $SO(3)_\mu=SO(3)$ if $\mu=0$. We look here for a number $N$ of vortices such that
$\dim\Fix(C_i,\p_N)=6\ (8\ \mbox{if}\ \mu=0)$ in order to obtain $\dim M_\mu=2$.
This case will be studied at the end of Section \ref{RPOnoident}.


\subsection{Relative periodic orbits formed of identical vortices}
\label{RPOident}

In the case of $N$ identical vortices, the subgroups $K\subset\Ker\chi$ lie in $SO(3)\times S_N$,
their elements are rotations coupled with cycles of $S_N$.
Moreover the knowledge of the spatial symmetries of an isotropy subgroup permits to fully describe
the subgroup, that is $K\simeq\pi(K)$ where $\pi:G\to O(3)$ is the Cartesian projection.

The Table \ref{lmrtabso3rpo} lists these subgroups:
the first column lists the groups using the usual
Sch\"{o}nflies notation; the second column gives the labels we use to
identify the different types of point-orbit of the action of that
group on $S^2$, the third column gives the isotropy subgroup for the
action at a point in that point-orbit and the fourth column the number
of points in the point-orbit;
the fifth column gives the dimension of $\Fix (K,\p_{|\oo|})$
(we give the number of connected components in parenthesis if it differs from one);
and the last column gives the dimension of $N_G(K)^o$.

A description of the different types of point-orbits is given in Appendix for each groups.
A configuration $x\in\p$ is labeled by $\pi(G_x)(k_1\oo_1,k_2\oo_2,\dots)$
where $\oo_j$ are the point-orbits of $G_x$ in $x$,
and $k_j$ number the occurrence of $\oo_j$ in $x$.
\begin{table}[tb]   
$$\begin{array}{|c|c|c|c|c|c|} \hline
 & & & & &  \\[-6pt]
K\simeq\pi(K) & \oo & I & |\oo| & \dim\Fix K & \dim N_G(K)^o \\[4pt] \hline
\hline \CC_n & R & 1   &n     &  2 &  1 \\
    & p & \CC_n &1    &  0(2)  &        \\
\hline \DD_n & R & 1   &2n  &  2  & 0 \\
    & r & \CC_2 &n   &  0(2)   &  \\
    & p & \CC_n &2   &  0  &  \\
\hline \Tet & R & 1    &12   & 2  & 0 \\
     & e & \CC_2  &6    &  0  &  \\
     & v & \CC_3  &4    &  0(2)  &  \\
\hline \Oct & R & 1    &24     &  2 & 0 \\
     & e & \CC_2  &12   &  0  &  \\
     & f & \CC_3  &8    &  0 &  \\
     & v & \CC_4  &6    &  0  &  \\
\hline \Icos & R & 1    &60    &  2   & 0 \\
      & e & \CC_2  &30  &  0 &  \\
      & f & \CC_3  &20  &  0 &  \\
      & v & \CC_5  &12  &  0 &  \\
\hline
\end{array}$$
\caption{Finite isotropy subgroups of $G=O(3)\times S_N$ formed of symplectic elements.}
\label{lmrtabso3rpo}
\end{table}

\bigskip

\noindent\textbf{Equilibria and Periodic orbits.}

\noindent A critical point of the restriction of $H$ to $\Fix K$ is a critical point of $H$ by the Principle
of Symmetric Criticality \cite{P79} since $H$ is $K$-invariant, and $K$ compact.
In particular, if $\dim\Fix K=0$, then points of $\Fix K$ are equilibrium points for the Hamiltonian system.
From Table \ref{lmrtabso3rpo}, configurations
$$C_n(p),D_n(r),\T(e),\T(v),\oct(e),\dots,\oct(v),\I(e),\dots$$
are therefore equilibrium configurations. In the same manner, we can show that
configurations $D_n(r,p),\T(e,v),\oct(e,f),\dots$ are also equilibria.
These results were already known \cite{LMR00} but show interest of Table \ref{lmrtabso3rpo}.

\bigskip

\noindent From the discussion of Section \ref{methode}, we expect periodic orbits in two-dimensional $\Fix K$.
Let $C_N$ act on $\p_{N+k_p}$ with $k_p=0,1$ or $2$, then $\Fix C_N$ is two-dimensional
and is formed of configurations
composed by a regular ring of $N$ identical vortices (a $N$-ring) together with $k_p$ polar vortices
of arbitrary vorticities. These configurations are \eqs from \cite{LMR00},
the fixed points manifold is therefore exclusively composed of \eqs (which are particular type of
periodic orbits since the dynamic orbit is closed).

Consider now $D_N$ acting on $\p_{2N+k_p}$ with $k_p=0,2$. The manifold $\Fix D_N$ is
two-dimensional, and composed of periodic orbits according to the study of \cite{ST}.
Their study shows also that this type of periodic orbits occurs for point vortex systems on a
cylinder or an ellipsoid of revolution.

Concerning polyhedra groups, the fixed points manifold is two-dimensional if $K=\T$ (resp.
$\oct$ and $\I$) acts on $\p_{12}$ (resp. $\p_{24}$ and $\p_{60}$), and is
formed of \ops \cite{ST}.
The \ops are composed of respectively $\T(R)$, $\oct(R)$ and
$\I(R)$ configurations which correspond to:
a regular tetrahedron (resp. octahedron and icosahedron)
for which each vertex is split in a $3$-ring (resp. $4$-ring and
$5$-ring), the vortices being on the edges starting from the vertex.
See \cite{ST} for a figure about the tetrahedral periodic orbit.

\bigskip

\noindent\textbf{Relative periodic orbits.}

\noindent From the discussion of Section \ref{methode}, $C_N$ is the only group that can be used in order to apply
the method (since $C_{h},C_{nh},C_i$ do not lie in $SO(3)$).

To obtain $\dim\Fix C_N=4$, we need to consider $C_N$ acting on $\p_{2N+k_p}$ with $k_p=0,1$ or $2$
(according to Table \ref{lmrtabso3rpo}). In that case, the fixed points manifold is formed of
$C_{Nv}(R_1,R_2,k_p p)$ configurations, that is composed of two
$N$-rings with $k_p$ polar vortices of arbitrary vorticities; and is indeed four-dimensional.
Consider the case without polar vortices, \ie $k_p=0$.
Assume without loss of generality that the vorticities of the two $N$-rings are
respectively $+1$ and $\lambda$. The fixed points manifold is locally parametrized by
$(\theta_1,\phi_1,\theta_{N+1},\phi_{N+1})$ where the vortices $x_1$ and $x_{N+1}$ belong
respectively to rings of vorticity $+1$ and $\lambda$.
We have then $\J_L=\tilde\J=N(\cos\theta_1+\lambda\cos\theta_{N+1})$,
the manifold $M_\mu$ is two-dimensional
 and is locally parametrized by $(\theta_1,\phi_{N+1})$, the two other variables satisfying $$
\phi_1=0,\ \theta_{N+1}=\arccos \left( \frac{\mu}{\lambda N} - \frac{\cos\theta_1}{\lambda}
\right). $$
A phase portrait in $M_\mu$ is given in Figure \ref{phaseCN}, it is obtained numerically
calculating (with \textsc{Maple}) the level sets of the function $\tilde H_\mu$.
\bfig \bc
\includegraphics[width=7cm,angle=0]{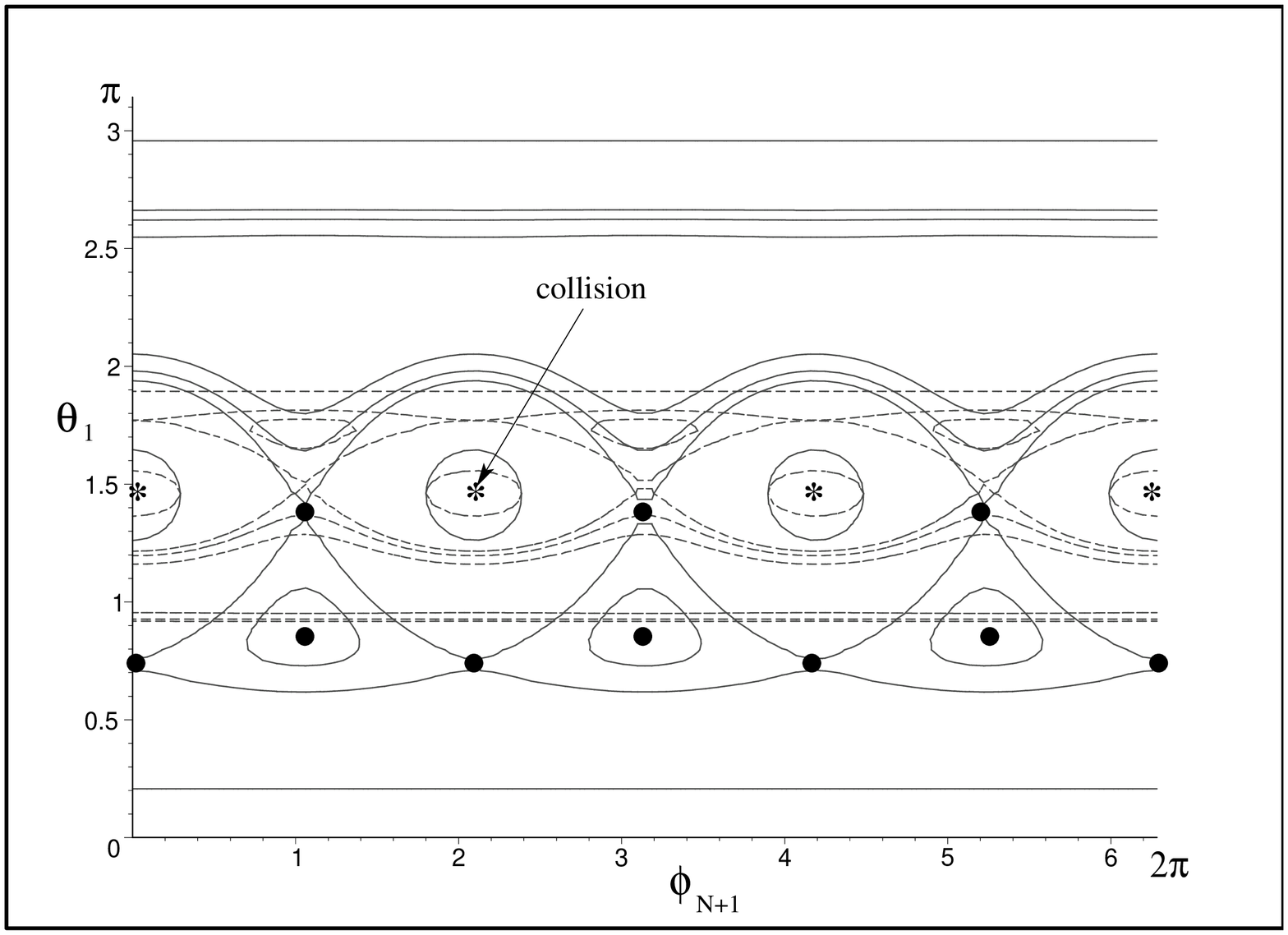}
\caption{Configurations  $C_{Nv}(R_1,R_2)$. Phase portrait in $M_\mu$ for $N=3$, $\lambda=2$
and $\mu=1$ (dots represent the trajectories for $\theta_{N+1}$).} \label{phaseCN}
\ec \efig


This phase portrait is composed of nine \eqs (three stable and six unstable), of RPOs, and of a
heteroclinic orbit connecting the six unstable \eqsb. There are three different families of RPOs:
those surrounding a stable \eqb, those surrounding a collision point, and those
which correspond to the -- more or less -- horizontal lines of the portrait.
For this last family of RPOs, the dynamics is the one of two uncoupled  $C_{Nv}(R)$ \ria with
angular velocities $\dot\phi_1$ and $\dot\phi_1+2\pi/T$ where $T$ is the period of
the RPO. For other values of $N,\lambda,\mu$, or if we add one or two polar vortices
of arbitrary vorticities ($k_p=1$ or $2$), we obtain RPOs too.\\

\rmk
It is important to note that the method lead to the existence of \ops and RPOs for a
\emph{arbitrary} number of vortices:
the phase portrait for other values of $N$ is similar to that for $N=3$.
This statement will go through for the phase portraits involving $C_N$ (or $C_{Nv}, C_{Nh}$) in
the next sections.

\subsection{Relative periodic orbits formed of vortices $(+)$ and $(-)$}
\label{RPOnoident}

\noindent We consider in that section $N$ vortices of vorticity $+1$ and
$N$ vortices of vorticity $-1$.
The interest is that $\Ker\chi$ is bigger in that case.
Indeed a reflection can here lie in $\Ker\chi$ (if it is coupled with the permutation $\tau$ which
exchange the vortices $(+1)$ with the vortices $(-1)$), contrary to the case of $N$ identical
vortices.

Due to the particular form of the phase space, we adapt the labels of the configurations.
Indeed, $K\simeq\pi(K)$ is no longer valid.
Let $K$ be an isotropy subgroup of $G=O(3)\times S_N \times S_N\rtimes \Z_2 \lbrack \tau \rbrack$,
the homomorphism $\varepsilon$ of Section \ref{pointsphere} takes the following form:
\begin{eqnarray*}
\varepsilon : {G} & \to & \lbrace -1,+1 \rbrace\\ (A,\sigma,\sigma^\prime,\tau^k) & \mapsto &
(-1)^k.
\end{eqnarray*}
We can then define $K_\omega=K\cap\Ker\chi$ and $K_\pm=K\cap\Ker\varepsilon$: $K_\omega$ is the
group of the symplectic elements of $K$, while $K_\pm$ is the group of the elements of $K$
preserving the sign of the vorticities. Note that $K_\omega$ and $K_\pm$ are subgroubs of index
one or two in $K$. We can remark also that to describe $K$ it is sufficient to know $\pi(K)$ and
$\pi(K_\pm)$ where $\pi$ is the Cartesian projection $\pi:G\to O(3)$. Thus we label the group $K$
by the pair $(\pi(K),\pi(K_\pm))$. Note that if $\pi(K)\subset SO(3)$, then $K_\pm=K_\omega$.
Moreover $K\subset\Ker\chi$ if and only if $K=K_\omega$.

The Table \ref{tableriasymp} lists the possible subgroups of $G$ which lie in $\Ker\chi$.
The way to label the orbits $\oo$ changes a little bit due again to the particular phase space.
For groups containing $C_n$, the orbits $\oo$ can be:
a ring of identical vortices ($R$),
a ring of alterned vortices ($\hat R$),
a semi-regular gon of $2n$ alterned vortices ($\hat R_s$),
two poles of opposite vorticities ($2p$).
Here \pg alterned vortices\pd means that the vortices $(+1)$ and the vortices $(-1)$
are alternatively placed on the orbit.
For polyhedral groups, we keep the labels of Table \ref{o(3)table} (see Appendix) when the vortices
are identical on the orbit, while
we add a hat on these labels when the vortices are alterned on the orbit.
For the group $C_i$, the orbit is a pair of antipodal vortices of opposite vorticities ($R$),
and for $C_h$, the orbit is a pair of vortices with opposite latitudes and vorticities,
but with same longitude ($R$).
\begin{table}[tb]
\begin{center}
\begin{tabular}{|c|c|c|c|c|c|c|}
\hline  $K\simeq (\pi(K),\pi(K_\pm))$ & $\pi(K_\omega)$ & $\oo$ & $I$ & $|\oo|$ & $\dim\Fix K$ & $\dim N_G(K)^o$\\ \hline
$(C_{n},C_{n})$ & $C_{n}$ & $2R$ & 1 & 2n & 4 & 1 \\
                &         & $2p$ & $C_{n}$ & 2 & 0(2) &  \\
\hline
$(C_{h},1)$ & $C_{h}$ & $nR$ & 1 & 2n & 2n & 1 \\
\hline
$(C_{nh},C_n)$ & $C_{nh}$ & $2R$ & 1 & 2n & 2 & 1 \\
  &  & $2p$ & $C_n$ & 2 & 0(2) &  \\
\hline
$(C_{i},1)$ & $C_i$ & $nR$ & 1 & 2n & 2n & 3 \\
\hline
$(C_{nv},C_n)$ & $C_{nv}$ & $\hat R_s$ & 1 & 2n & 2 & $0$ \\
$(D_{nd},D_{n})$ & $D_{nd}$ & $2\hat R_s$ & 1 & 4n & 2 & $0$ \\
$(D_{nh},D_{n})$ & $D_{nh}$ & $2R_s$ & 1 & 4n & 2 & $0$ \\
$(\T,\T)$ & $\T$ & $v+v'$ & $C_{3}$ & 8 & 0(2) &  $0$\\
$(\T_h,\T)$ & $\T_h$ & $\hat v$ & $C_{3}$ & 8 & 0(2) &  $0$\\
$(\T_h,\T)$ & $\T_h$ & $\hat R$ & 1 & 24 & 2 &  $0$\\
$(\T_d,\T)$ & $\T_d$ & $\hat R$ & 1 & 24 & 2 & $0$\\
$(\oct_h,\oct)$ & $\oct_h$ & $\hat R$ & 1 & 48 & 2 &  $0$\\
$(\I_h,\I)$ & $\I_h$ & $\hat R$ & 1 & 120 & 2 &  $0$\\
$(\T,\T)$ & $\T$ & $2R$ & 1 & 48 & 4 &  $0$\\
$(\oct,\oct)$ & $\oct$ & $2R$ & 1 & 96 & 4 & $0$\\
$(\I,\I)$ & $\I$ & $2R$ & 1 & 240 & 4 &  $0$\\
\hline
\end{tabular}
\caption{Finite subgroups of $G=O(3)\times S_N \times
S_N\rtimes \Z_2 \lbrack \tau\rbrack$ formed of symplectic elements.}
\label{tableriasymp}
\end{center}
\end{table}

Following the previous section, we label a configuration $x\in\p$ by
$$(\pi(K),\pi(K_\pm))(k_1\oo_1,k_2\oo_2,\dots)$$
where $K=G_x\simeq(\pi(K),\pi(K_\pm))$, $\oo_j$ are the point-orbits of $G_x$ in $x$,
and $k_j$ number the occurrence of $\oo_j$ in $x$.
Two different labels can lead to the same arrangement
(this because a subgroup is not necessarily an isotropy subgroup).

\bigskip

\noindent\textbf{Equilibria and Periodic orbits.}

\noindent As in Section \ref{RPOident}, points of $\Fix K$ with $\dim\Fix K=0$ are equilibria.
It follows that configurations $(C_{n},C_{n})(2p),(C_{nh},C_{n})(2p)$ and
$(\T,\T)(v+v^\prime),(\T_h,\T)(\hat v)$ are equilibria.
The two first correspond to the well-known equilibrium configuration formed of two polar vortices
(here with opposite vorticities).
The two second correspond to another known equilibrium \cite{LP02}:
a regular tetrahedron formed of four vortices $(+1)$ together with its tetrahedron dual
formed of vortices $(-1)$. That arrangement can also be seen as a cube formed of alterned vortices
($\hat v$).

\bigskip

We expect here also periodic orbits in two-dimensional $\Fix K$.
The configurations $(C_{h},1)(R)$ satisfy $\dim\Fix G_x=2$, they correspond to
pair of vortices with opposite latitude and vorticities, and with same longitude.
These configurations are well-known \ria \cite{KN98}.
The same hold for $(C_{nh},C_n)(2R)$: these configurations are \ria \cite{LP02},
they are formed of a $n$-ring of vortices $(+1)$ together with a $n$-ring of vortices $(-1)$, the
two rings being at opposite latitudes and \pg aligned\pd (in phase).
For that two cases, the two-dimensional fixed point manifold $\Fix K$ is exclusively composed of
relative equilibria.

Configurations $(C_{Nv},C_N)(\hat R_s)$ correspond to $2N$ vortices placed at the vertices of
a semi-regular polygon at fixed latitude, the vortices being alterned.
The fixed points manifold $\Fix K$ is two-dimensional and we can easily obtain
the phase portrait of the Hamiltonian subsystem $(\Fix K,H|_{\Fix K})$.
The manifold $\Fix K$ is locally parametrized by coordinates of one vortex, say $(\theta_1,\phi_1)$.
One has $(\theta_1,\phi_1)\in (0,\pi)\times (-\pi/(2N),\pi/(2N))$ since collisions are excluded.
The phase portrait is composed of one-parameter family of \ops (Figures \ref{dance2} and \ref{dessdance}),
these \ops are the \emph{Dancing vortices} trajectories of \cite{To}.
\bfig \bc
\subfigure[]{\label{dance2}\includegraphics[width=6cm,angle=0]{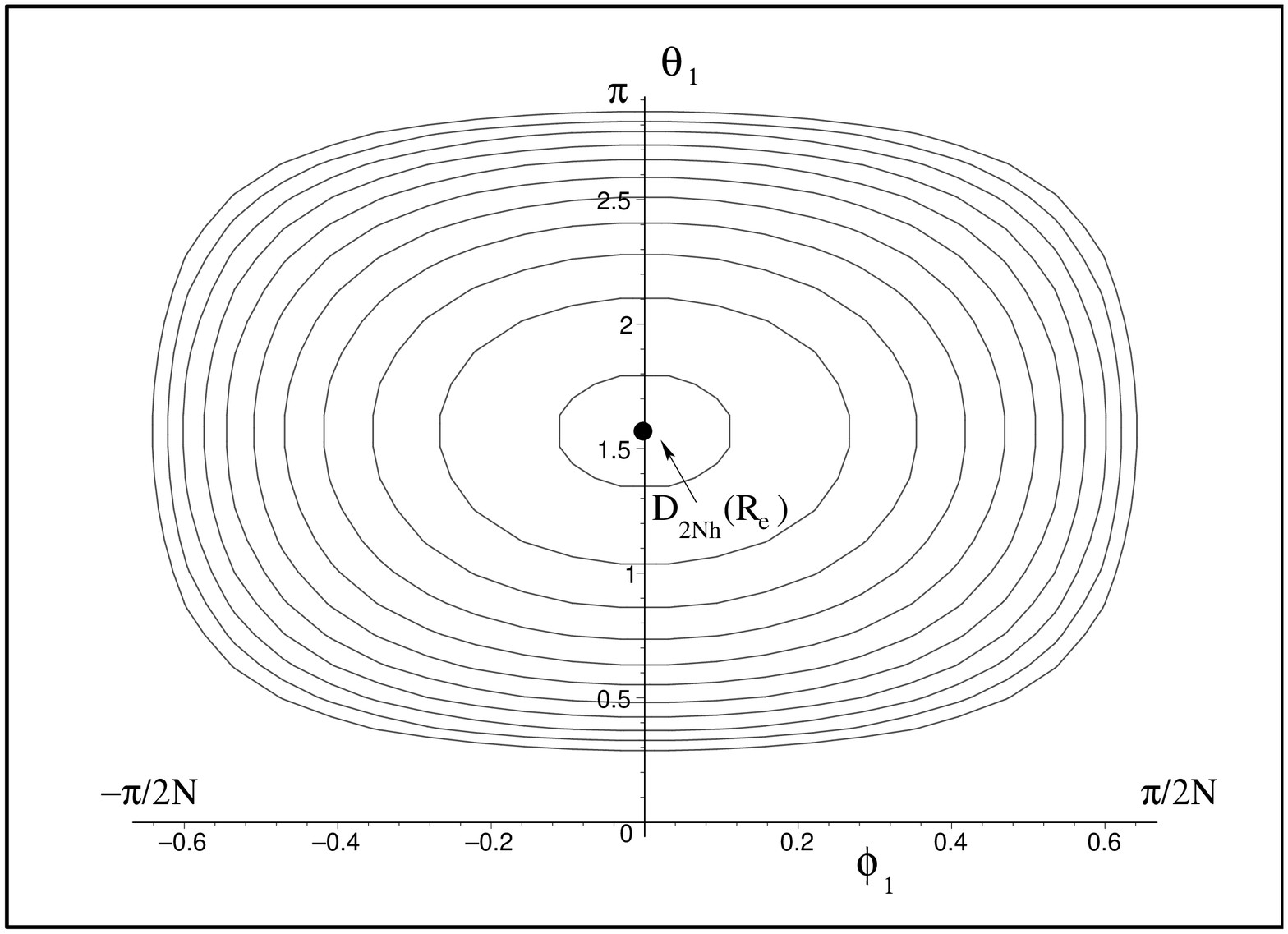}}
\hspace{0.6cm}
\subfigure[]{\label{dessdance}\includegraphics[width=4cm,angle=0]{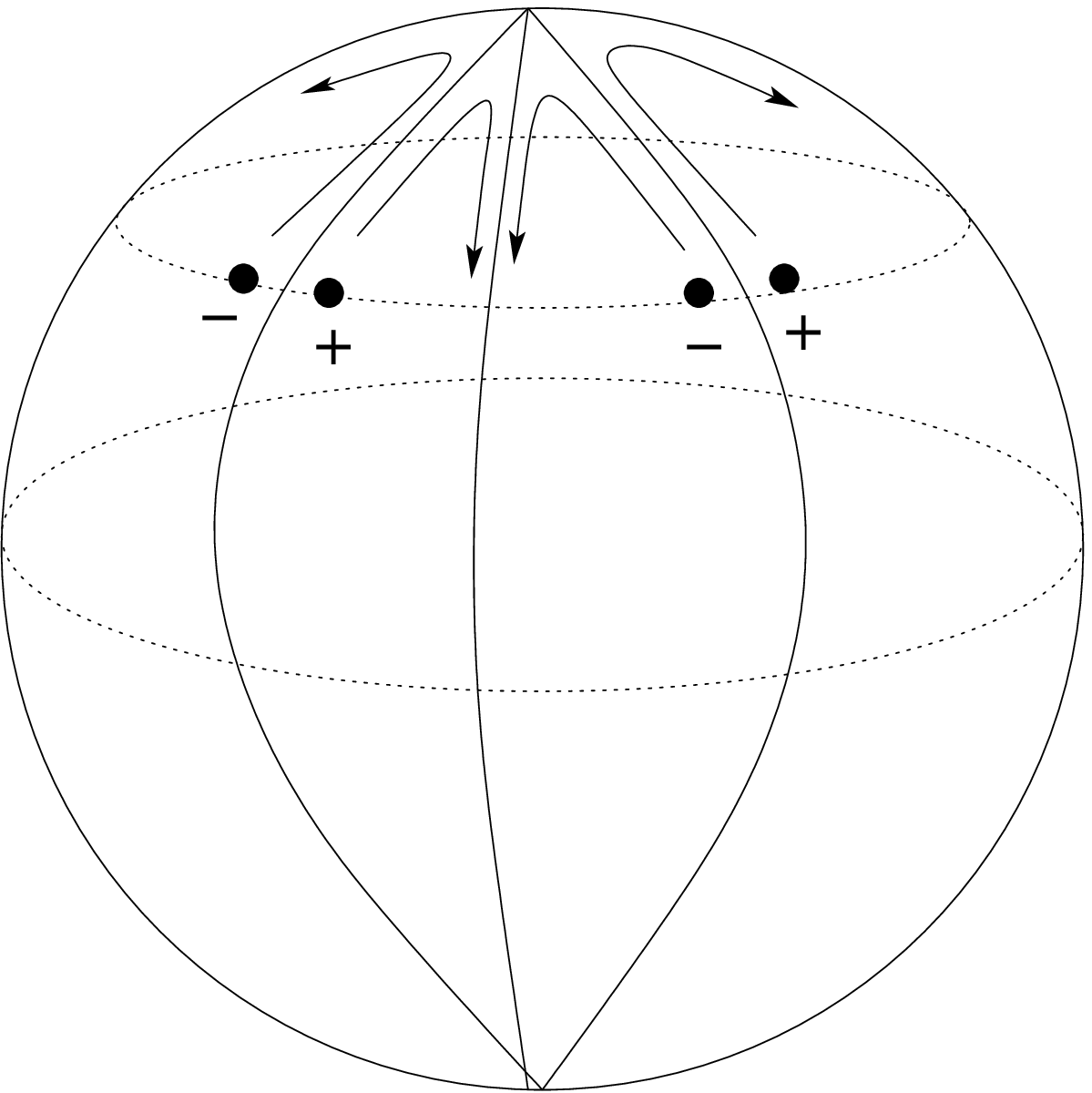}}

\caption{(a) Phase portrait for \emph{Dancing vortices} with $N=2$, the portrait is similar for $N>2$.
(b) Motion of \emph{Dancing vortices} for $N=4$, only the front hemisphere is represented.}
\ec \efig
These \ops surround the equilibrium point $D_{2Nh}(R_e)$
(which is an equatorial ring formed of $2N$ alterned vortices) which is Lyapunov stable (modulo $SO(3)$)
 if $N=2$, unstable if $N>2$ \cite{LP02}. The stability of these \ops has been calculated in \cite{LPTh}:
the stability of the $D_{2Nh}(R_e)$ equilibrium persists along the family of \opsb,
\ie the \ops are stable if $N=2$, unstable if $N>2$.

Consider now $4N$ vortex of vorticities $\lambda_j=+1,\lambda_{2N+j}=-1,j=1\dots2N$.
Start with configurations $(D_{Nd},D_{N})(2\hat R_s)$. The fixed points manifold $\Fix G_x$
is two-dimensional and
is locally parametrized by $(\theta_0,\phi_0)\in(0,\pi)\times(0,\pi/N)$ in such way that
$$\theta_{j}=\theta_{2N+j}=\theta_0,\theta_{N+j}=\theta_{3N+j}=\pi-\theta_0,j=\1N$$ $$
\phi_{j}=\frac{2\pi j}{N}-\phi_0, \phi_{N+j}=\frac{2\pi j}{N}+\frac{\pi}{N}+\phi_0,
\phi_{2N+j}=\frac{2\pi j}{N}+\phi_0, \phi_{3N+j}=\frac{2\pi j}{N}+\frac{\pi}{N}-\phi_0,j=\1N $$
for all configurations in $\Fix G_x$. The phase portrait is composed exclusively of \ops
surrounding the collision point $(\theta_0=\pi/2,\phi_0=\pi/(2N))$ (Figure \ref{conf6}).

In the same manner, configurations $(D_{nh},D_{n})(2R_s)$ lie on periodic orbits,
$\Fix G_x$ is locally parametrized by
$(\theta_0,\phi_0)\in(0,\pi/2)\times(-\pi/(2N),\pi/(2N))$ in such a way that
$$\theta_{j}=\theta_{2N+j}=\theta_0,\theta_{N+j}=\theta_{3N+j}=\pi-\theta_0,j=\1N$$ $$
\phi_{j}=\frac{2\pi j}{N}+\phi_0, \phi_{N+j}=\frac{2\pi j}{N}+\frac{\pi}{N}-\phi_0,
\phi_{2N+j}=\frac{2\pi j}{N}+\frac{\pi}{N}-\phi_0, \phi_{3N+j}=\frac{2\pi j}{N}+\phi_0,j=\1N.$$
The phase portrait is composed of \ops surrounding an equilibrium point (Figure \ref{conf7}).
\bfig \bc
\subfigure[]{\label{conf6}\includegraphics[width=5.6cm,angle=0]{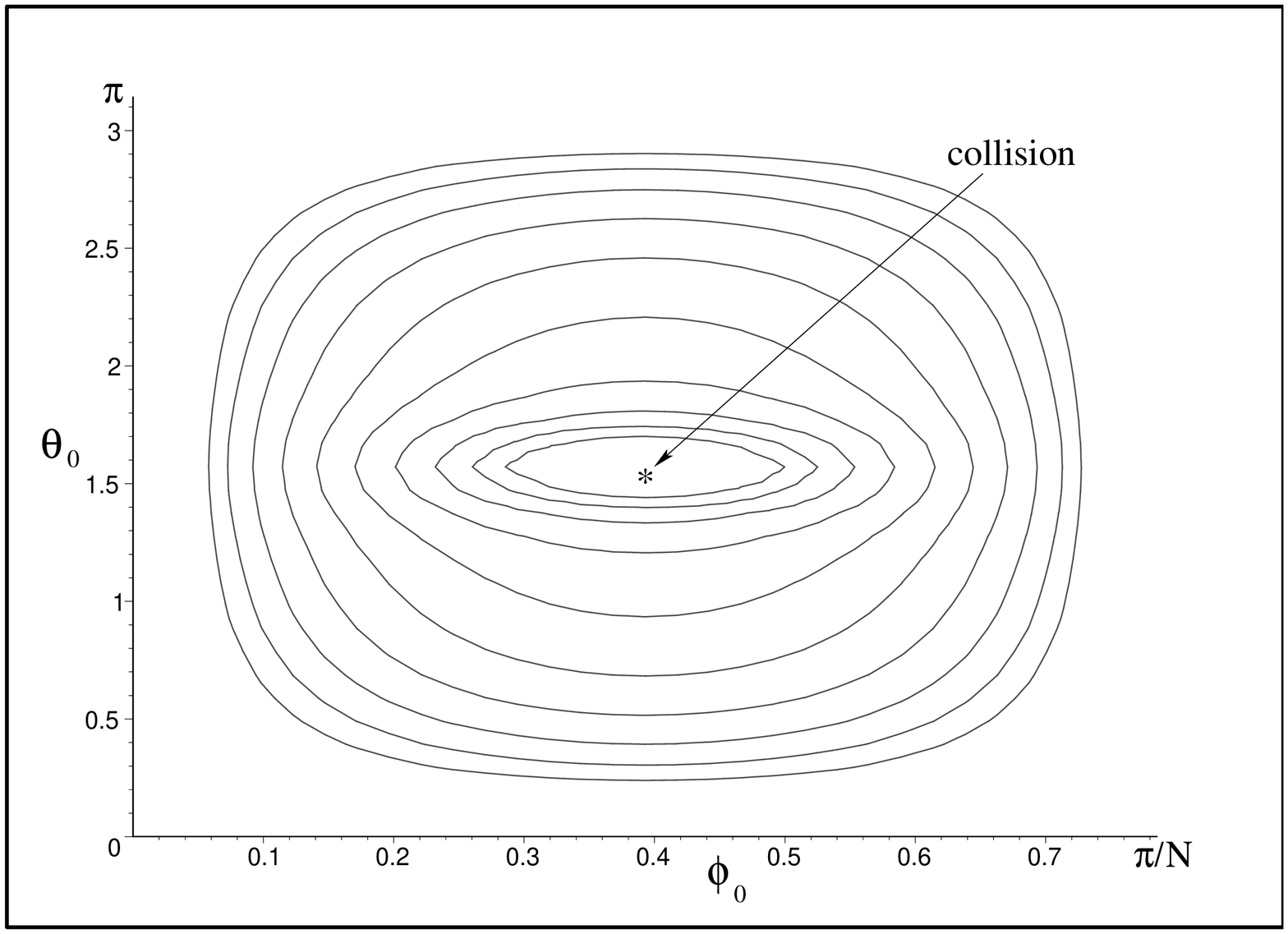}}
\hspace{0.1cm}
\subfigure[]{\label{conf7}\includegraphics[width=5.6cm,angle=0]{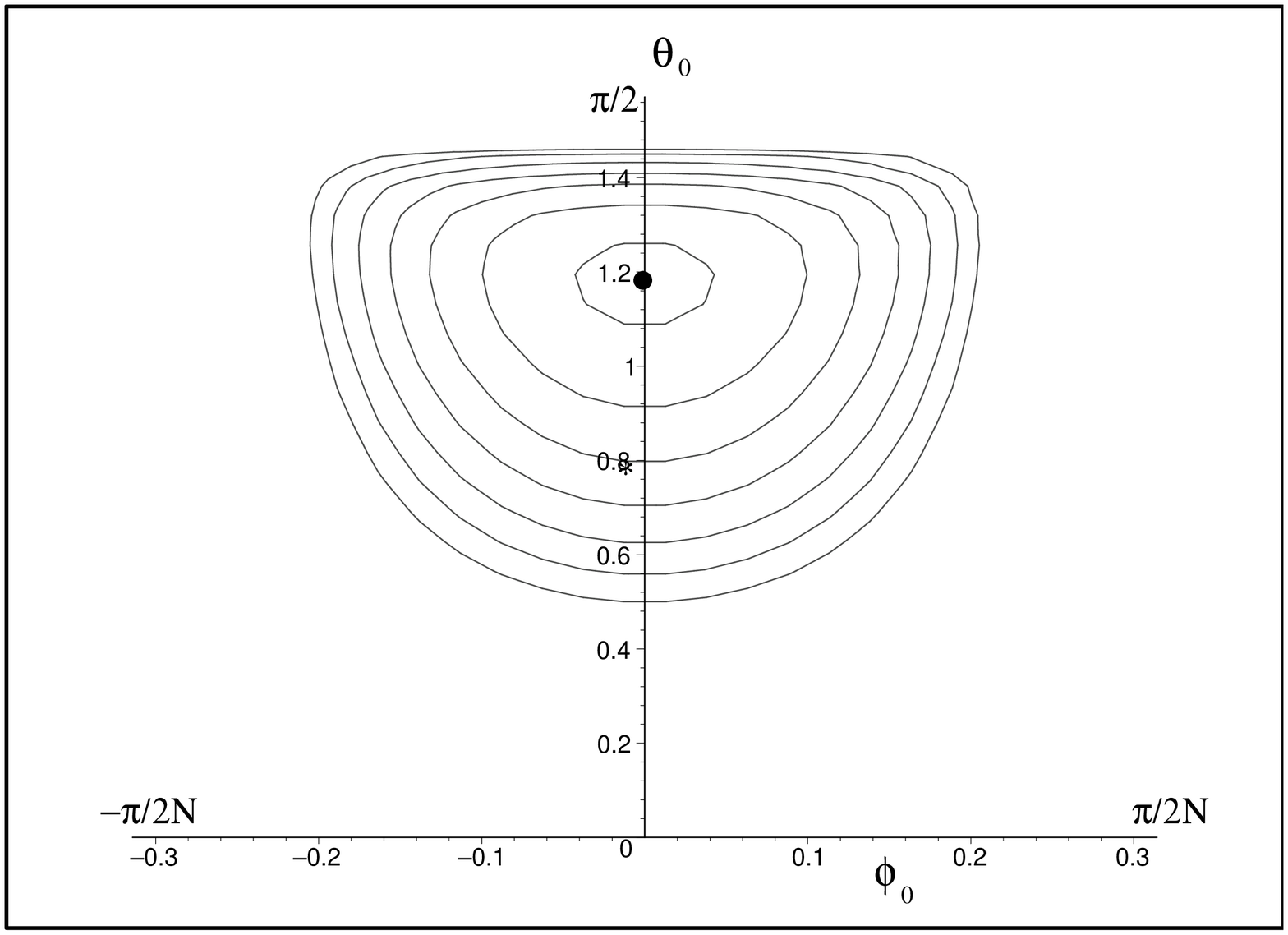}}
\caption{(a) Phase portrait for configurations $(D_{Nd},D_{N})(2\hat R_s)$ with $N=4$.
(b) Phase portrait for configurations $(D_{Nh},D_{N})(2R_s)$ with $N=5$.
The portraits are similar -- respectively -- for other values of $N$.}
\ec \efig
We have in particular exhibited an equilibrium composed of two alterned $2N$-rings
(i.e. of type $2\hat R$) with opposite latitudes.

To end with periodic orbits, consider the polyhedral groups.
Configurations $(\T_h,\T)(\hat R)$ correspond to a regular tetrahedron of
vortices $(+1)$ and its dual formed of vortices $(-1)$ for which each vertex has been split in a $3$-ring,
the three vortices being on three edges starting from the vertex.
For configurations $(\T_d,\T)(\hat R)$ (resp. $(\oct_h,\oct)(\hat R)$ and $(\I_h,\I)(\hat R)$),
each vertex of the regular tetrahedron (resp. regular octahedron and regular icosahedron) is split in
a semi-regular $2\times3$ (resp. $2\times4$ and $2\times5$)-ring of alterned vortices
(see Figure \ref{splitvertex}).
The fixed points manifolds for these different cases are two-dimensional. Moreover one can easily verify
that $\tilde H$ is proper and non constant, hence these manifolds contain \ops
from Proposition \ref{proprenoncste}.
\bfig \bc
\includegraphics[width=3cm,angle=0]{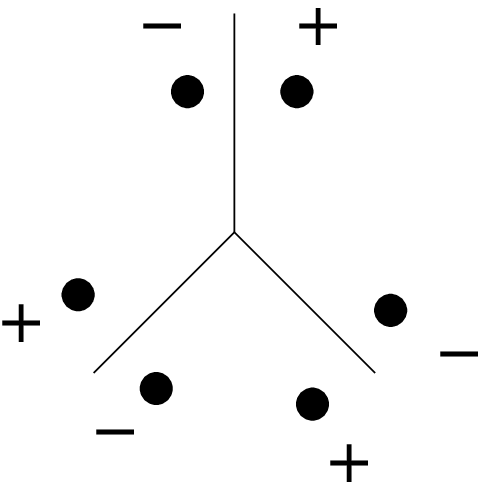}
\caption{Splitting of a vertex of the regular tetrahedron in a semi-regular $2\times3$-ring  of
alterned vortices (view from top).} \label{splitvertex} \ec \efig

\bigskip

\noindent\textbf{Relative periodic orbits.}
We dispose of the four groups $C_h,C_n,C_{nh},C_i$ to perform the method of Section \ref{methode}.

\bigskip

\noindent\textbf{$C_h$ symmetry.}
Let $K=(C_h,1)$ act on $\p_4$ the phase space of two vortex $(+1)$ with two vortex $(-1)$.
The manifold $\Fix K$ is four-dimensional, and locally parametrized by the coordinates $(\theta_1,\phi_1,\theta_2,\phi_2)$
of the two vortices $(+1)$. One has $\J_L=\tilde\J=2(\cos\theta_1+\cos\theta_2)$, then we can
parametrize $M_\mu$ by $(\theta_1,\phi_1)$, the other variables of $\p$ satisfying $$ \phi_2=0,\
\theta_{2}=\arccos \left( \frac{\mu}{2} - \cos\theta_1 \right),\ \phi_3=\phi_1 ,\
\theta_3=\pi-\theta_1,\ \phi_4=0,\ \theta_4=\pi-\theta_2.$$
A phase portrait is given in Figure \ref{chrpofig} for $\mu=0.8$,
we distinguish two families of relative periodic orbits:
those surrounding the stable \eqb, and those surrounding the poles.
When one of the vortices $(+1)$ moves along a trajectory above the $1-3$ collision line
($\theta_1=\pi/2$), the other moves along a trajectory below the $2-4$ collision line.
The stable \eq occurs for $\theta_1=\theta_2$ and $\phi_1=\pi$, which corresponds
to the $C_{2v}(R_m,R'_m)$ \eq of \cite{LP02}.
\bfig \bc
\includegraphics[width=7cm,angle=0]{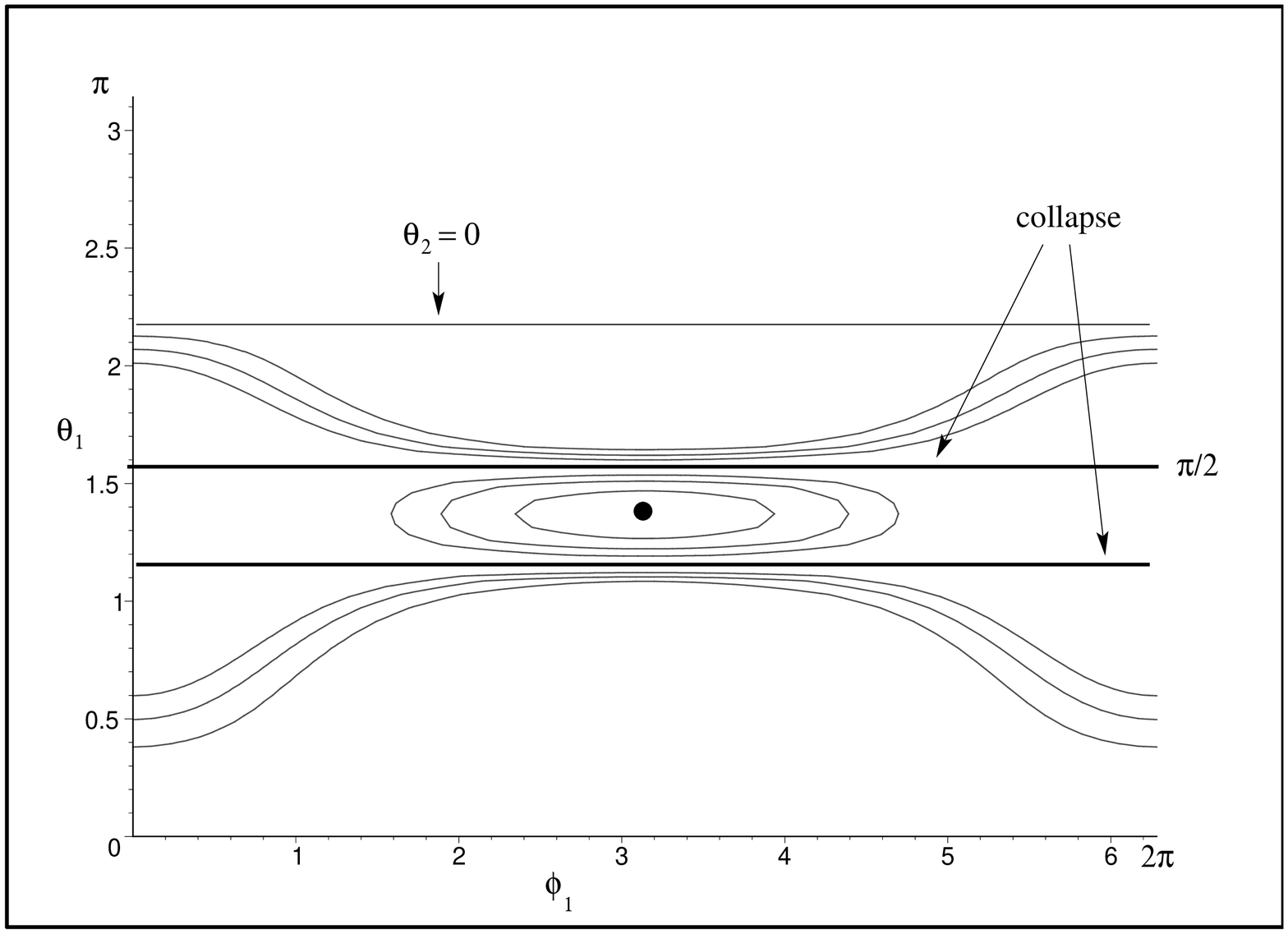}
\caption{Phase portrait for configurations $(C_h,1)(2R)$ with $\mu=0.8$.}
\label{chrpofig} \ec\efig

\bigskip

\noindent\textbf{$C_n$ symmetry.}
Let $K=(C_N,C_N)$ act on the phase space of $2N+2$ vortices.
The manifold $\Fix K$ is four-dimensional, and corresponds to configurations $(C_N,C_N)(2R,2p)$
formed of two polar vortices of opposite vorticities together with two $N$-rings of arbitrary latitudes and
of vorticity respectively $+1$ and $-1$.

Let $(\theta_1,\phi_1)$ (resp. $(\theta_1^\prime,\phi_1^\prime)$) be the coordinates of a
vortex on the ring of vorticity $+1$ (resp. $-1$), the manifold $M=\Fix K$ is locally parametrized by
$(\theta_1,\phi_1,\theta_1^\prime,\phi_1^\prime)$. One has
$\J_L=\tilde\J=N(\cos\theta_1-\cos\theta_1^\prime)+2\lambda$ where $\lambda$ is the vorticity of the
North pole. We have $\dim M_\mu=2$ and $M_\mu$ is locally parametrized by $(\theta_1,\phi_1^\prime)$,
the two other variables of $M_\mu$ satisfying $$\phi_1=0,\
\theta_1^\prime=\arccos \left( \frac{2\lambda-\mu}{N} + \cos\theta_1 \right). $$
Take $\mu=2\lambda$, hence $\theta_1=\theta_1^\prime$, the two $N$-rings form a
semi-regular $2N$-ring as in $(C_{Nv},C_N)(\hat R_s)$ configurations.
The phase portrait in $M_\mu$ is similar to Figure \ref{dance2} whatever the value of $\lambda$, the trajectories
of this portrait correspond therefore to relative periodic orbits.
Actually, the value of $\lambda$ does not affect the level curves of $\tilde H_\mu$ but only on their
\pg energies\pd ($\tilde H_{\mu,\lambda}=\tilde H_\mu-\lambda^2\ln2$).
To summarize, we have shown that the \emph{Dancing vortices} \ops
become RPOs when two polar vortices of opposite vorticities are added.

\bigskip

\rmk
This last statement is actually intuitive:
adding two polar vortices of opposite vorticities to the $D_{2Nh}(R_e)$ equilibrium configuration
(an equatorial ring of $2N$ identical vortices) makes the whole arrangement a relative
equilibrium;
hence the idea that adding two polar vortices of opposite vorticities changes a periodic orbit into
a RPO.


\bigskip

For $\mu\neq2\lambda$, we obtain more exotic RPOs. Fix $N=3$ and $\lambda=+1$ for simplicity.

For $\mu=4$, the phase portrait is the one of a planar pendulum (Figure \ref{phasesN3mu4}),
we have therefore shown the existence of \eqs (three stable and three unstable), RPOs, and of a
homoclinic cycle connecting the three unstable \eqs
(which shape actually a single configuration).
\bfig \bc
\includegraphics[width=7cm,angle=0]{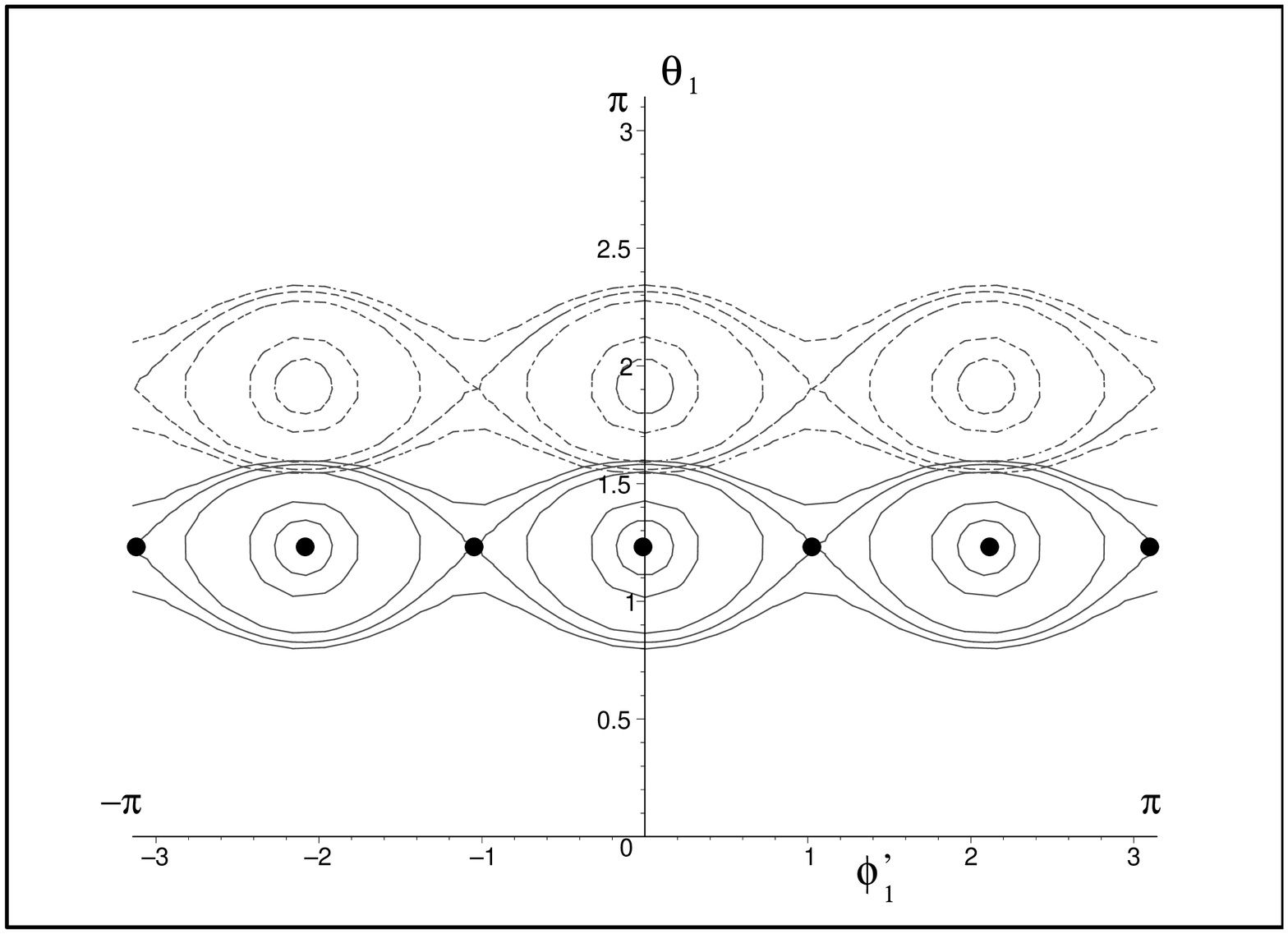}
\caption{Configurations $(C_N,C_N)(2R,2p)$.
Phase portrait in $M_\mu$ for $N=3$, $\lambda=+1$ and $\mu=4$
(dots represent the trajectories for $\theta_1^\prime$).} \label{phasesN3mu4} \ec
\efig
This portrait is not specific to the value $\mu=4$, it is similar for $\mu=3.1$ for example.
Take now $\mu=2.5$, we obtain the phase portrait of the Figure \ref{muzerohuit},
hence we show the existence of 12 \eqs (six stable and six unstable), of a homoclinic cycle
connecting the six unstable \riab, and of two types of RPOs. The first type of RPO
corresponds to those surrounding the stable \rumb, the second type corresponds to the almost
horizontal lines of the phase portrait (in a frame rotating with the ring $(+)$, the ring
$(-)$ rotates around itself). A bifurcation must \emph{a priori} occur between $\mu=2.5$ and $\mu=3.1$.
Actually, a sub-critical pitchfork bifurcation occurs between these two values when $\mu$ decreases.
\bfig \bc
\includegraphics[width=7cm,angle=0]{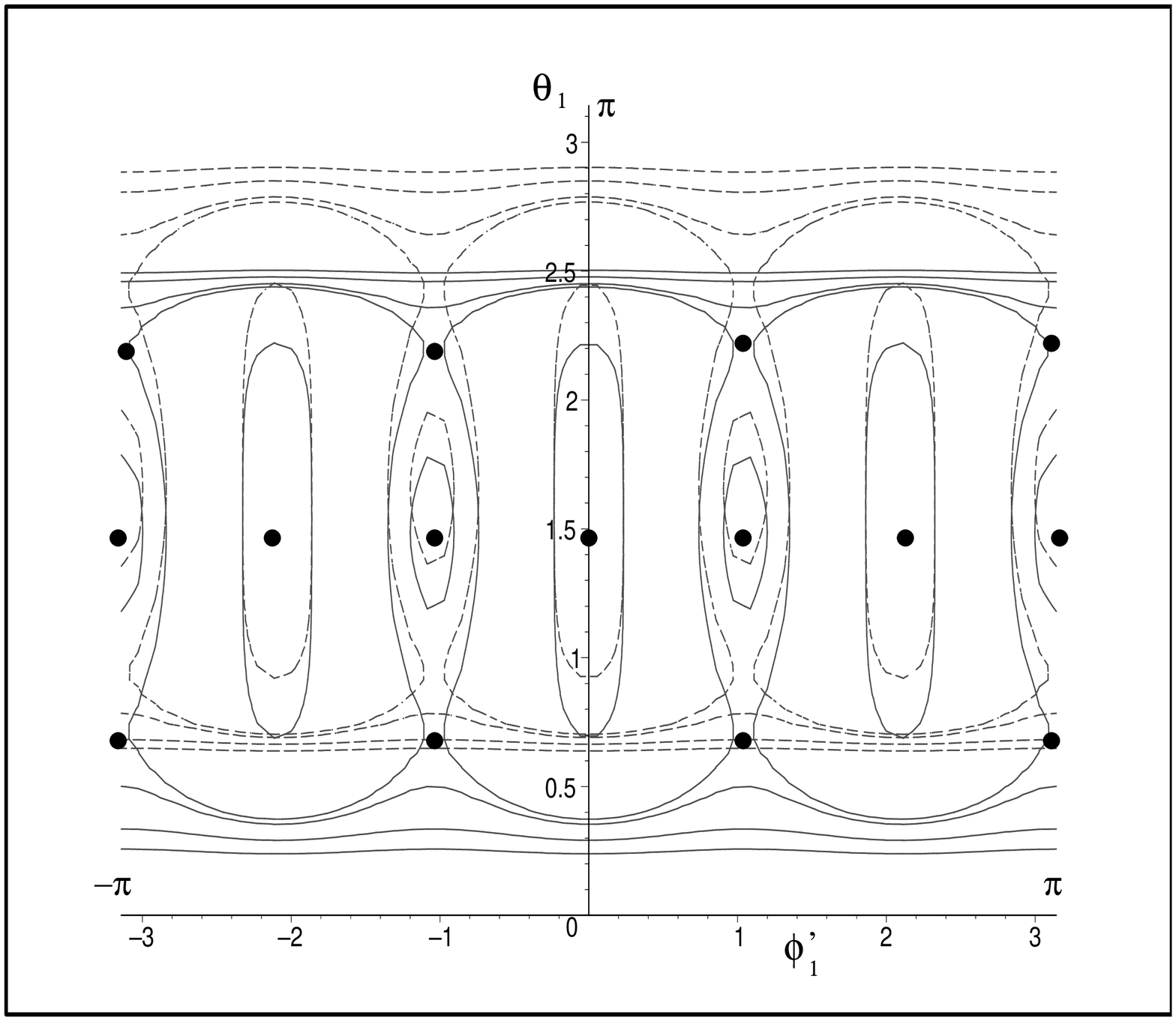}
\caption{Configurations $(C_N,C_N)(2R,2p)$.
Phase portrait in $M_\mu$ for $N=3$, $\lambda=+1$ and $\mu=2.5$
(dots represent the trajectories for $\theta_1^\prime$).} \label{muzerohuit} \ec
\efig
Calculus of phase portraits show that these phenomena persist for other values of $N$ and $\lambda$.\\

\bigskip

\noindent\textbf{$C_{nh}$ symmetry.}
Let $K=(C_{Nh},C_N)$ act on the phase space of $4N+2$ vortices.
The manifold $\Fix K$ is four-dimensional, and corresponds to configurations $(C_{Nh},C_N)(4R,2p)$
formed of two polar vortices of opposite vorticities, together with two $N$-rings of vorticity $(+1)$
and of arbitrary latitudes, together with two $N$-rings of vorticity $(-1)$ and of latitudes opposite
to those of the rings of vorticity $(+1)$.

The manifold $M=\Fix K$ is locally parametrized by $\theta_1,\phi_1,\theta_{2N+1},\phi_{2N+1}$
where the vortices $x_1$ and $x_{2N+1}$ do not belong to the same orbit of $C_{Nh}$,
and $\lambda_1=-\lambda_{2N+1}=1$. With this setting, one has
$\J_L=\tilde\J=2N(\cos\theta_1-\cos\theta_{2N+1})+2\lambda$,
where $\lambda$ is the vorticity of the North pole.
The manifold $M_\mu$ is two-dimensional and locally parametrized by $(\theta_0,\phi_0)$
where $\theta_0=\theta_1$ and $\phi_0=\phi_{2N+1}$, the two other variables
of $M_\mu$ satisfying $$ \phi_1=0,\ \theta_{2N+1}=\arccos \left( \frac{2\lambda-\mu}{2N} +
\cos\theta_0 \right).$$
We can assume in addition $\phi_0\!\in(0,2\pi/N)$ without loss of generality.

Take $\mu=2\lambda$, hence $\theta_{2N+1}=\theta_1$, the four $N$-rings form
two semi-regular $2N$-rings ($2\hat R_s$) with opposite latitudes.
The phase portrait in $M_\mu$ is given Figure \ref{conf7rpomu2},
we obtain RPOs and two \eqsb. These \eqs correspond to two $2N$-rings $\hat R$ with opposite latitudes
(this for a specific latitude).
Here again the phase portrait is independent of the value for $\lambda$ since
$\tilde H_{\mu,\lambda}=\tilde H_\mu-\lambda^2\ln2$.
\bfig \bc
\includegraphics[width=7cm,angle=0]{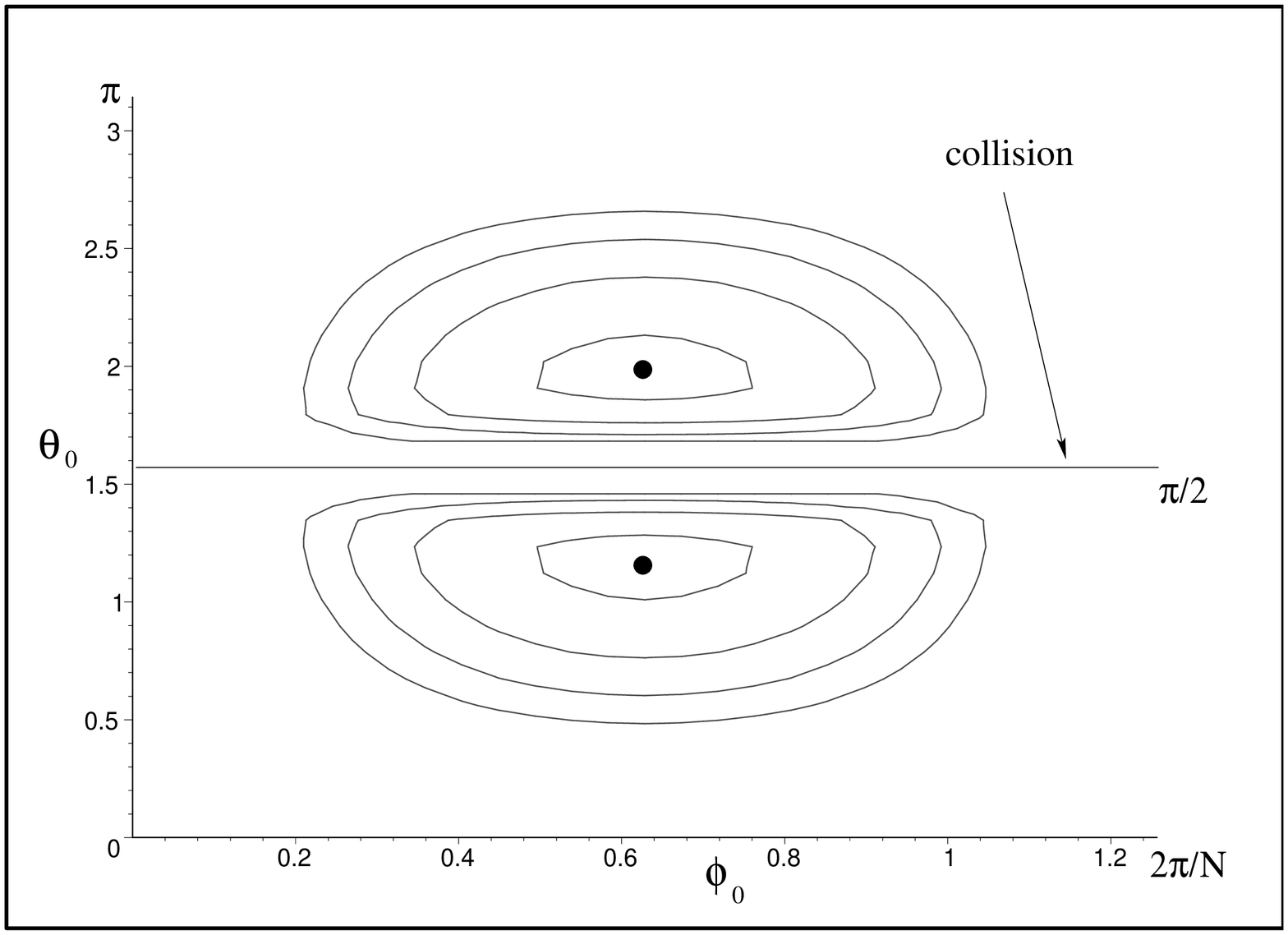}
\caption{Configurations $(C_{Nh},C_N)(4R,2p)$.
Phase portrait for $\mu=2\lambda$. The portrait is for $N=5$.}
\label{conf7rpomu2} \ec \efig

Consider now cases for which $\mu\neq2\lambda$. Fix $N=3$ and $\lambda=+1$.

For $\mu=6$, the phase portrait (Figure \ref{conf7rpomu6}) is composed of three unstable relative equilibria,
of RPOs, and of a homoclinic cycle connecting the three unstable \eqs
(which shape actually a single configuration). Here again RPOs are of two types: those surrounding
a collision point, and those for which $\theta_0$ is almost constant.
\bfig \bc
\includegraphics[width=7cm,angle=0]{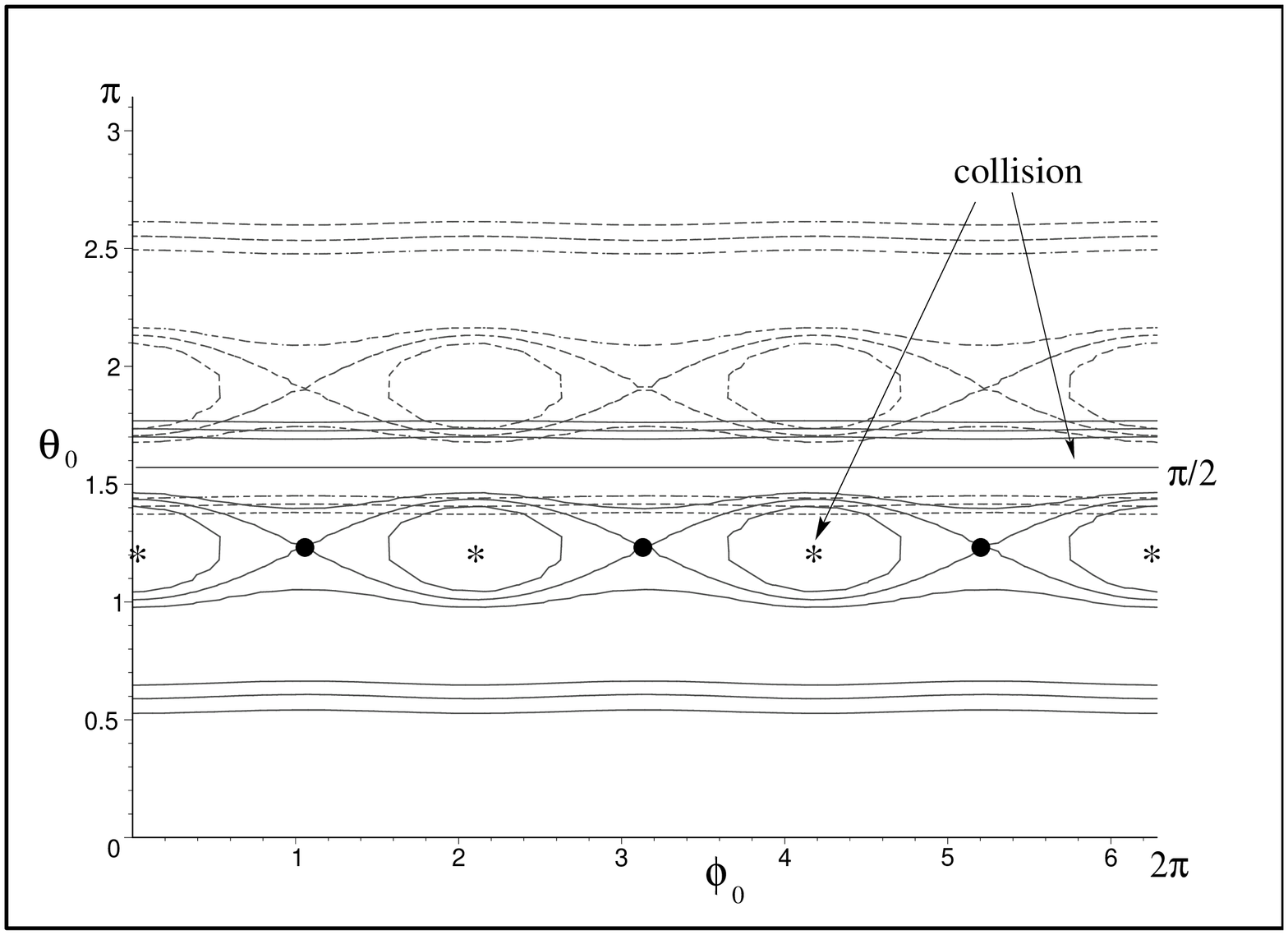}
\caption{Configurations $(C_{Nh},C_N)(4R,2p)$.
Phase portrait in $M_\mu$ for $N=3$, $\lambda=+1$ and $\mu=6$
(dots represent the trajectories for $\theta_{2N+1}$).} \label{conf7rpomu6} \ec
\efig

For $\mu=2.5$, the phase portrait is slightly different (Figure \ref{conf7rpomu25}):
it is composed of twelve \eqs (six stable and six unstable), of RPOs, and of two homoclinic cycles
each connecting three unstable \eqsb. These homoclinic cycles are remarkable:
the cycle is formed of a first cycle connecting the three \eqs and of three homoclinic orbit
starting from each of the three \eqsb.
However, it is not clear what type of bifurcation occurs between $\mu=2.5$ and $\mu=6$,
\pg saddle-centre\pd bifurcations play probably a role.
\bfig \bc
\includegraphics[width=7cm,angle=0]{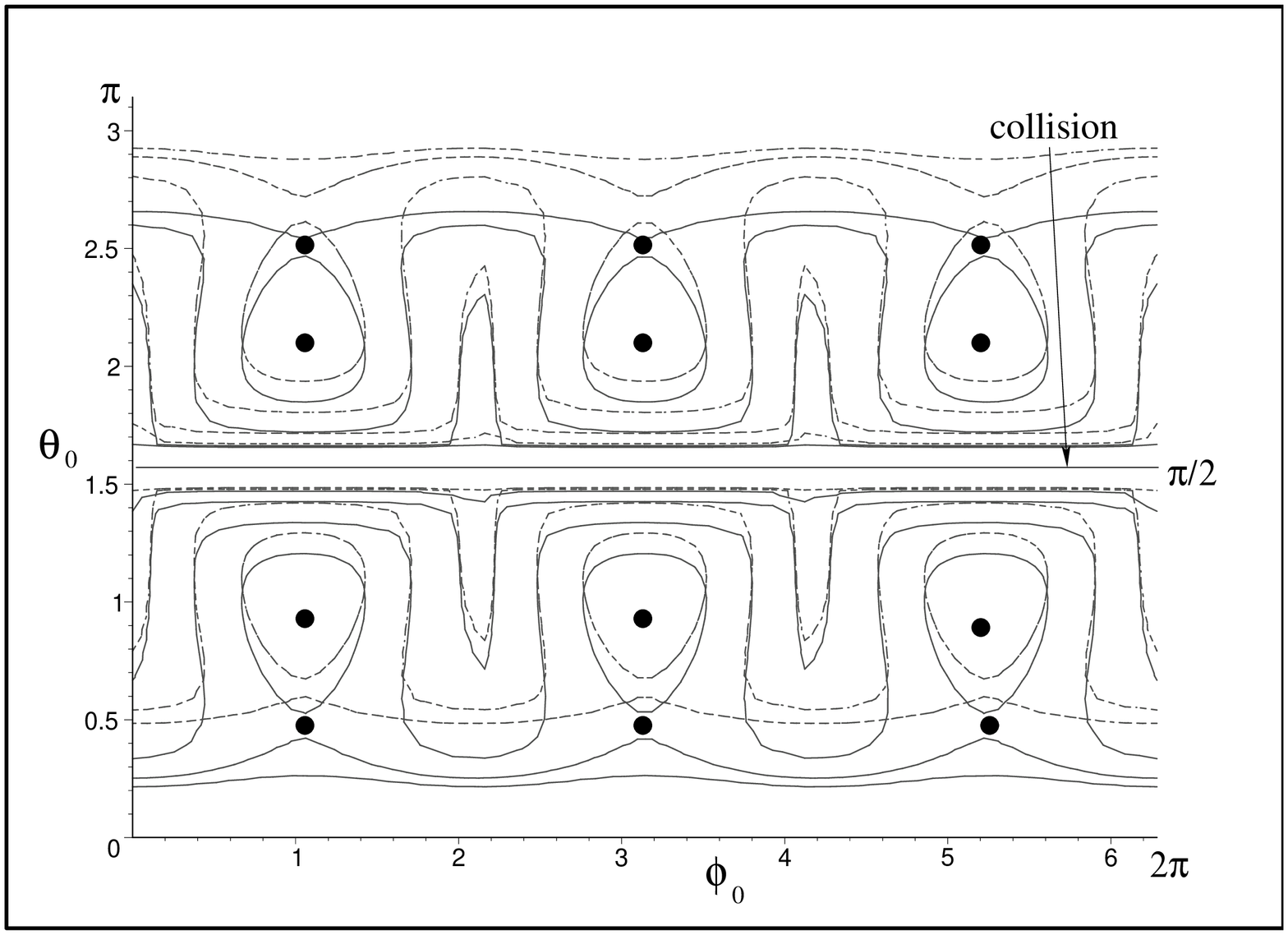}
\caption{Configurations $(C_{Nh},C_N)(4R,2p)$.
Phase portrait in $M_\mu$ for $N=3$, $\lambda=+1$ and $\mu=2.5$
(dots represent the trajectories for $\theta_{2N+1}$).} \label{conf7rpomu25} \ec
\efig

We show calculating other phase portraits that these phenomena persist for other values of $N$ and $\lambda$.\\

\bigskip

\noindent\textbf{$C_i$ symmetry.}
To end, we study the case $K=(C_i,1)$.
Let act $(C_i,1)$ on the phase space of six vortices
(three vortices $(+1)$ and three vortices $(-1)$). We have $\Fix K\simeq S^2\times S^2\times
S^2\setminus \Delta$ where the three spheres correspond to the three vortex $x_1,x_2,x_3$ of
vorticities $+1$ and $\Delta$ is the set of all possible collisions in $\Fix C_i$, hence
$\dim\Fix K=6$. The symplectic form in $\Fix C_i$ is given by $$
\tilde{\omega}=2\sin\theta_1\ d\theta_1\land d\phi_1+2\sin\theta_2\ d\theta_2\land d\phi_2
+2\sin\theta_3\ d\theta_3\land d\phi_3 $$ and $L=N_G(K)^o=SO(3)$ acts diagonally on $\Fix K$,
hence we get
$$
\begin{array}{rcl}
\J_L:\Fix K & \to & \liel^*\simeq\R^3\\ (x_1,x_2,x_3) & \mapsto & 2(x_1+x_2+x_3)
\end{array}
$$ (embedding $S^2$ in $\R^3$).

We have $\J_L=\tilde\J$ and so $M_\mu=\tilde\J^{-1}(\mu)/SO(3)_\mu$.
Zero-momentum configurations of $\Fix K$ are equilateral triangles lying in a great circle
(this configuration in $\p$ is the $D_{6h}(R_e)$ equilibrium).

Assume $\mu\neq 0$. We have $SO(3)_\mu=SO(2)$, the manifold $M_\mu$ is therefore two-dimensional,
and is locally parametrized by the co-latitudes $(\theta_1,\theta_2)$ of vortices $x_1$ and $x_2$.
The other vortex coordinates are given by $$ \theta_3=\arccos\left(
\frac{\mu}{2}-\cos\theta_1-\cos\theta_2 \right),\ \theta_4=\pi-\theta_1 ,\ \theta_5=\pi-\theta_2,\
\theta_6=\pi-\theta_3 $$ $$ \phi_1=0,\ \phi_2=\arccos\left( -
\frac{\sin^2\theta_1+\sin^2\theta_2-s_3^2}{2\sin\theta_1\sin\theta_2} \right) $$ $$
\phi_3=-\arccos\left( - \frac{\sin^2\theta_1-\sin^2\theta_2+s_3^2}{2s_3\sin\theta_1} \right),\
\phi_4=\pi+\phi_1 ,\ \phi_5=\pi+\phi_2,\ \phi_6=\pi+\phi_3 $$ with $$ s_3=\sqrt{1-\left(
\frac{\mu}{2}-\cos\theta_1-\cos\theta_2 \right)^2}, $$
which permit to reconstruct trajectories in $M_\mu$ into the original phase space $\p$.
The phase portrait in $M_\mu$ for $\mu=1.5$ is given in Figure \ref{cirpo}, it is composed of
a family of RPOs surrounding the stable \eq $D_{3d}(R,R')$, and of a heteroclinic cycle connecting
the unstable \ria $D_{2h}(2R,2p)$. This heteroclinic cycle is formed of great circle
configurations; in a frame rotating with that great circle, the
vortices are following each other on the same trajectory.\\
\bfig \bc
\includegraphics[width=7cm,angle=0]{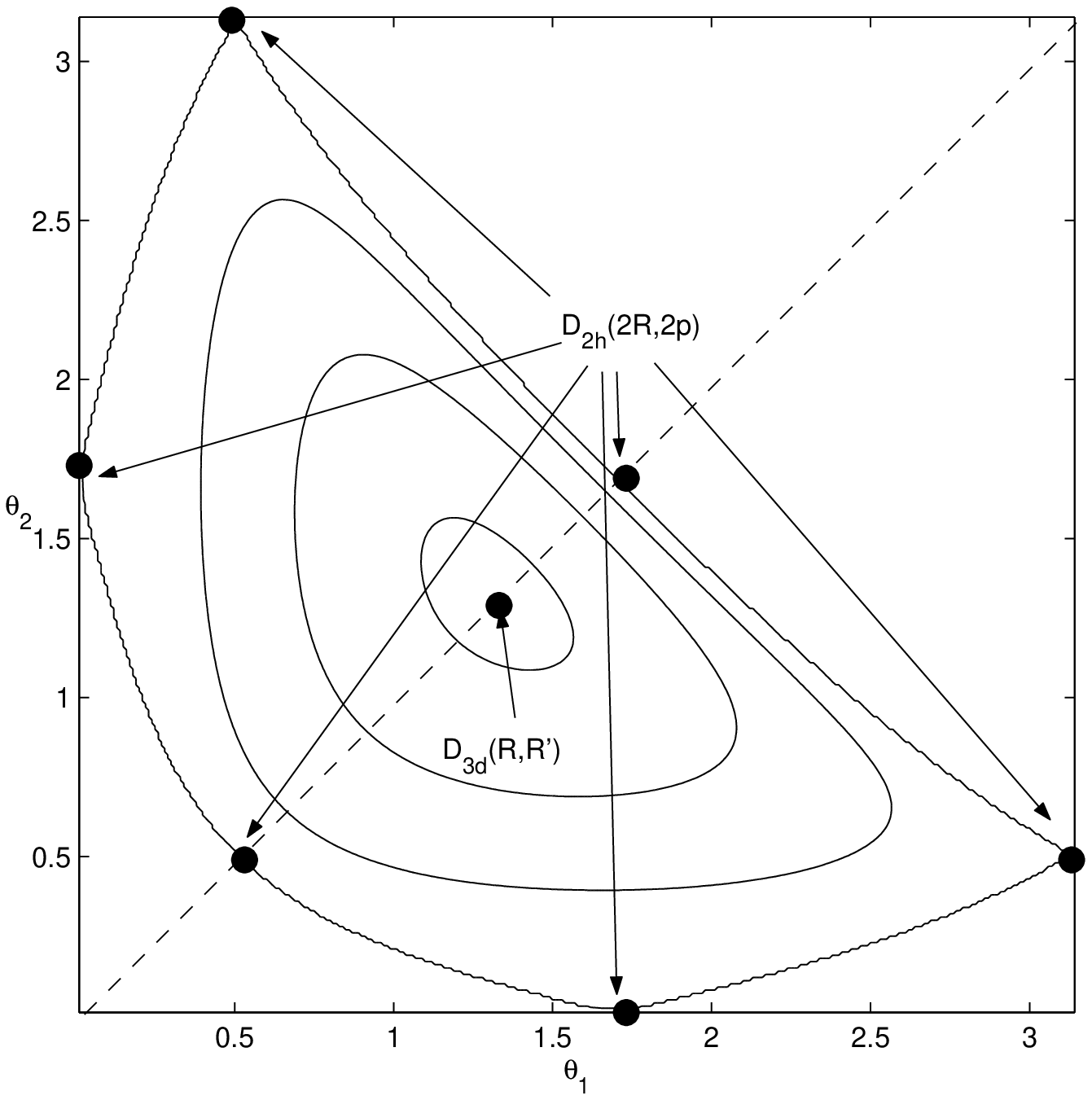}
\caption{Configurations $(C_i,1)(3R)$. Phase portrait in $M_\mu$ for $\mu=1.5$.}
\label{cirpo} \ec \efig

There is also a particular case for which we can obtain $\dim M_\mu=2$ with eight vortices.
Consider configurations $(C_i,1)(4R)$ with $\mu=0$.
We have $SO(3)_\mu=SO(3)$, thus $\dim M_0=2$.
After some straightforward calculus, we find that $M_0$ can be locally parametrized by the coordinates
$(\theta_3,\phi_3)\in (0,\pi)\times(\frac{\pi}{2},\frac{3\pi}{2})$ of $x_3$,
the other variables satisfying
$$
\theta_1=\phi_1=\phi_2=0,\
\theta_2=\arccos\left(
\frac{(1-\cos^2\theta_3)\cos^2\phi_3-(1+\cos\theta_3)^2}{(1-\cos^2\theta_3)\cos^2\phi_3+(1+\cos\theta_3)^2}
\right)\in\left[\frac{\pi}{2},\pi\right),
$$
$$x_4=-(x_1+x_2+x_3),\ x_5=-x_1,\ x_6=-x_2,\ x_7=-x_3,\ x_8=-x_4,$$
which permit to reconstruct trajectories in $M_0$ into the original phase space $\p$.
The phase portrait (Figure \ref{cirpo8}) is formed of RPOs surrounding a collision point.
\bfig \bc
\includegraphics[width=7cm,angle=0]{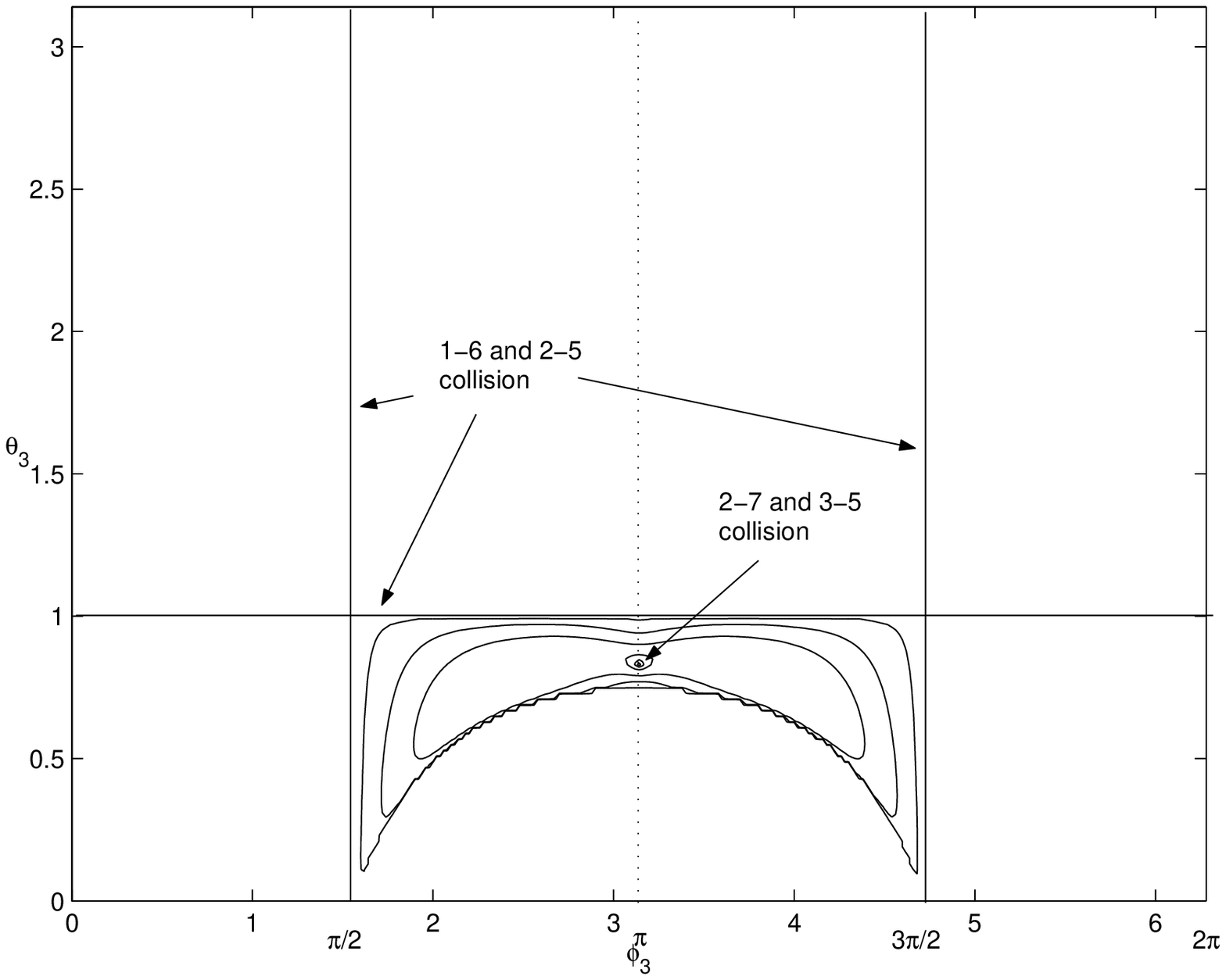}
\caption{Configurations $(C_i,1)(4R)$. Phase portrait in $M_{\mu=0}$.}
\label{cirpo8} \ec \efig

\rmk
The work of the former two sections can also be done to determine RPOs in systems of point
vortices on a surface with symmetry (ellispoid, hyperboloid, cylinder \dots):

$\bullet$ the $C_i$ symmetry study will work for surfaces with $\Z_2[-Id]$ symmetry,

$\bullet$ the $C_h$ symmetry study will work for surfaces with $\Z_2[z\mapsto -z]$ symmetry,

$\bullet$ the $C_n$ symmetry study will work for surfaces with a symmetry axis,

$\bullet$ the $C_{nh}$ symmetry study will work for surfaces with a symmetry axis together
          with a $\Z_2[z\mapsto -z]$ symmetry.

We expect the two-dimensional phase portraits on these surfaces to give RPOs and (relative)
heteroclinic cycles.

\section{Point vortices in the plane}
\label{rpoplane}

\subsection{Description of the Hamiltonian system}
\label{pointplane}

The equations of motion $N$ planar point vortices $z_1,\dots,z_N\in\C$ are \cite{Ar82,LR96}:
$$\dot{\overline{z_j}}=\frac{1}{2\pi i}\sum_{k=1,k\neq j}^{N} \frac{\lambda_k}{z_j-z_k}$$ where
$\lambda_j$ is the vorticity of the vortex $z_j$.
The Hamiltonian for this system
is $$ H=-\frac{1}{4\pi}\sum_{i<j} \lambda_i \lambda_j \ln\vert z_i-z_j \vert^2 $$
and the vector field $X_H$ is equivariant under the action of $SE(2)$ (which is not compact).
After identifying $SE(2)$ with $\C\rtimes SO(2)$ and so $\se^*$ with $\C\times\R$,
the \mom of the system is $$\J=\left( i\sum_{j=1}^{N}
\lambda_j z_j , \frac{1}{2}\sum_{j=1}^{N} \lambda_j \rho_j^2 \right)$$
where $z_j=\rho_j\exp (i\phi_j)$ for all $j=1,\dots,N$.

In order to apply the method of Section \ref{methode}, we are interested only in compact symmetry
groups, thus we forget translational symmetries: the vector field $X_H$ is $SO(2)$-equivariant
and the \mom due to that $SO(2)$-symmetry is $$\J=\frac{1}{2}\sum_{j=1}^{N} \lambda_j \rho_j^2$$
which is one conserved quantity.

Following Section \ref{pointsphere}, the Hamiltonian $H$ is invariant under $G=O(2)\times\hat S(\loup)$,
and the vector field $X_H$ is $\Ker(\chi)$-equivariant, where $S(\loup)$ and $\chi$ are defined as in
Section \ref{pointsphere}.
In the case of $N$ identical vortices, the Hamiltonian is $O(2)\times S_N$-invariant,
while $X_H$ is $SO(2)\times S_N$-equivariant.

\subsection{Application of the method}

Compare to the problem on the sphere, the planar problem is less rich since we dispose of only two
types of finite groups: $C_n$ and $D_n$.
From these, the only groups with a continuous normalizer are $C_n$ groups: $N_G(C_n)^o=SO(2)$.
Since $C_n\subset SO(2)$, we will not consider as in Section \ref{RPOnoident} the particular
phase space of $N$ vortices $(+1)$ with $N$ vortices $(-1)$ and its remarkable symmetries.

The point-orbits of $C_N$ can be a regular $N$-ring ($R$) or a central vortex ($p$).
Consider $N$ identical vortices, we have $\dim\Fix G_x=2$ if $x$ is a $C_N(R)$ configuration;
while $\dim\Fix G_x=0$ for a $C_N(p)$ configuration.

To obtain $\dim M_\mu=2$, we need to get $\dim\Fix K=4$ since $SO(2)_\mu=SO(2)$ ($SO(2)$ acts
trivially on $\sodeux^*\simeq\R$).
To this end, we consider $N$ vortices of vorticity $+1$ together with $N$ vortices of vorticity
$\lambda$ and a central vortex of vorticity $\lambda_c$. In that case $\Fix (C_N,\p_{2N+1})$
is four-dimensional and formed of $C_N(2R,p)$ configurations.
The fixed points manifold can be parametrized by coordinates of two vortices $z_1,z_{N+1}$ which belong to
different rings.
We have $\J_L=\tilde\J=N(\rho_1^2+\lambda\rho_{N+1}^2)/2$ and $\dim M_\mu=2$,
the manifold $M_\mu$ can be parametrized by $(\rho_1,\phi_1)$, the other variables of $M=\Fix K$
satisfying $\rho_{N+1}=\sqrt{(2\mu/N-\rho_1^2)/\lambda}$ and $\phi_{N+1}=0$ ($\phi_1$ is actually
the offset between the two rings).
Note that if $\lambda>0$, then the dynamics are necessarily bounded since
$\rho_{1}=\sqrt{2\mu/N-\lambda\rho_{N+1}^2}$.

Different phase portraits are given in Figure \ref{planarfig} for $N=3$ and $\mu=1$.
Actually, the value of the momentum $\mu$ does not affect the shape of the phase portrait, it just gives
a characteristic length scale.
Note also that, as in the spherical case,
the value of the vorticity $\lambda_c$ of the central vortex does not affect the level curves of
the phase portrait but only their energies.

For $\lambda=-5$, the phase portrait (Figure \ref{cnm5}) is composed of RPOs:
in a frame rotating with the ring of vorticity $(+1)$, the ring of vorticity
$\lambda$ rotates around itself.
For $\lambda=-1$, it is composed of RPOs together with unbounded trajectories (Figure \ref{cnm1}).
For $\lambda=+1$, the phase portrait (Figure \ref{cn1}) shows the existence of 9 \eqs
(three stable and six unstable), of a heteroclinic cycle
connecting the six unstable \riab, and of three types of RPOs. The first type of RPO
corresponds to those surrounding a stable \rumb, the second to those surrounding a collision
point, and the third type corresponds to the almost
horizontal lines of the phase portrait (for which the motion is as in Figure \ref{cnm5}).

\bfig \bc
\subfigure[$\lambda=-5$]{\label{cnm5}\includegraphics[width=5cm,angle=0]{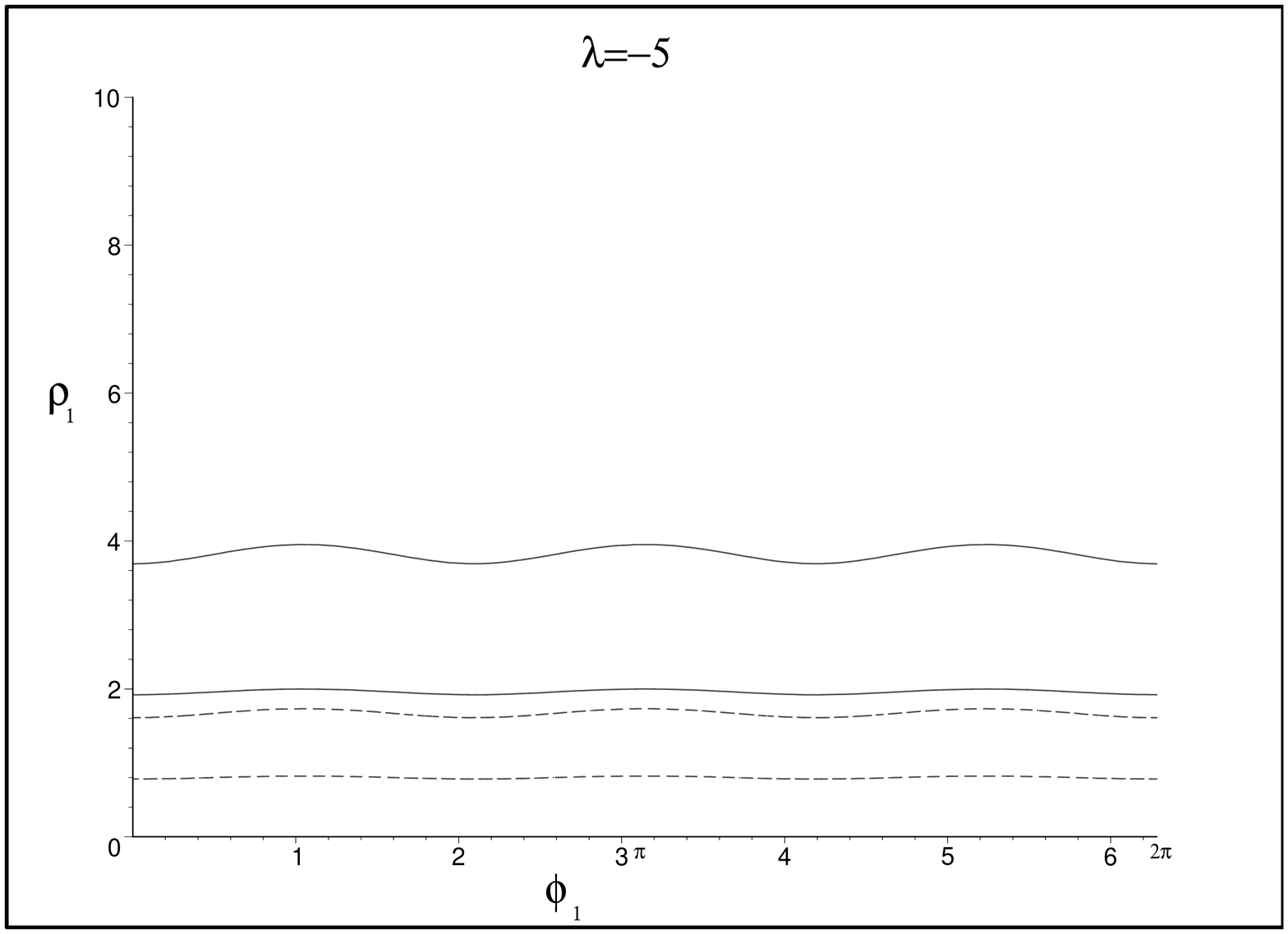}}
\hspace{0.6cm}
\subfigure[$\lambda=-1$]{\label{cnm1}\includegraphics[width=5cm,angle=0]{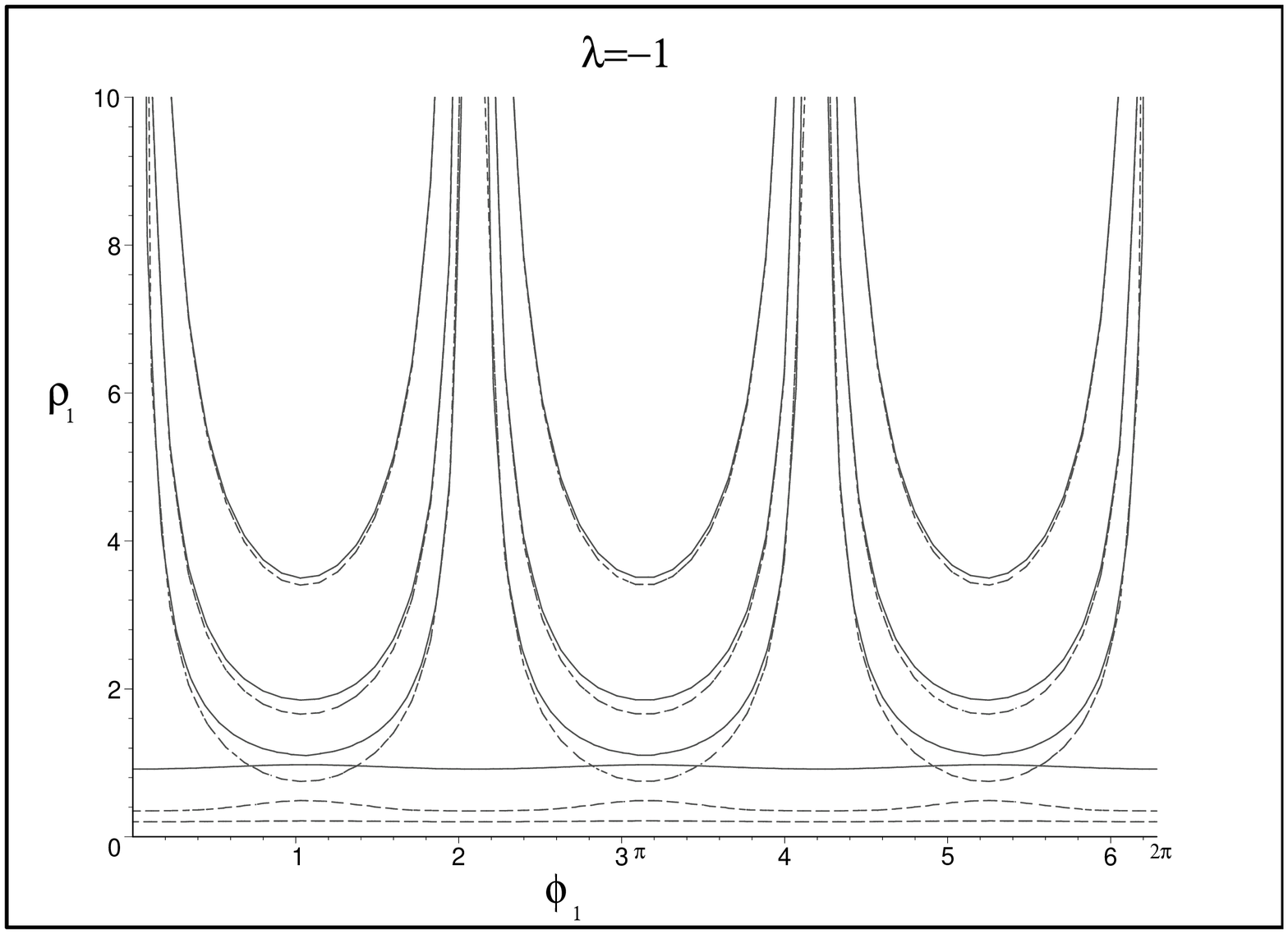}}

\vspace{0.6cm}

\subfigure[$\lambda=+1$]{\label{cn1}\includegraphics[width=7cm,angle=0]{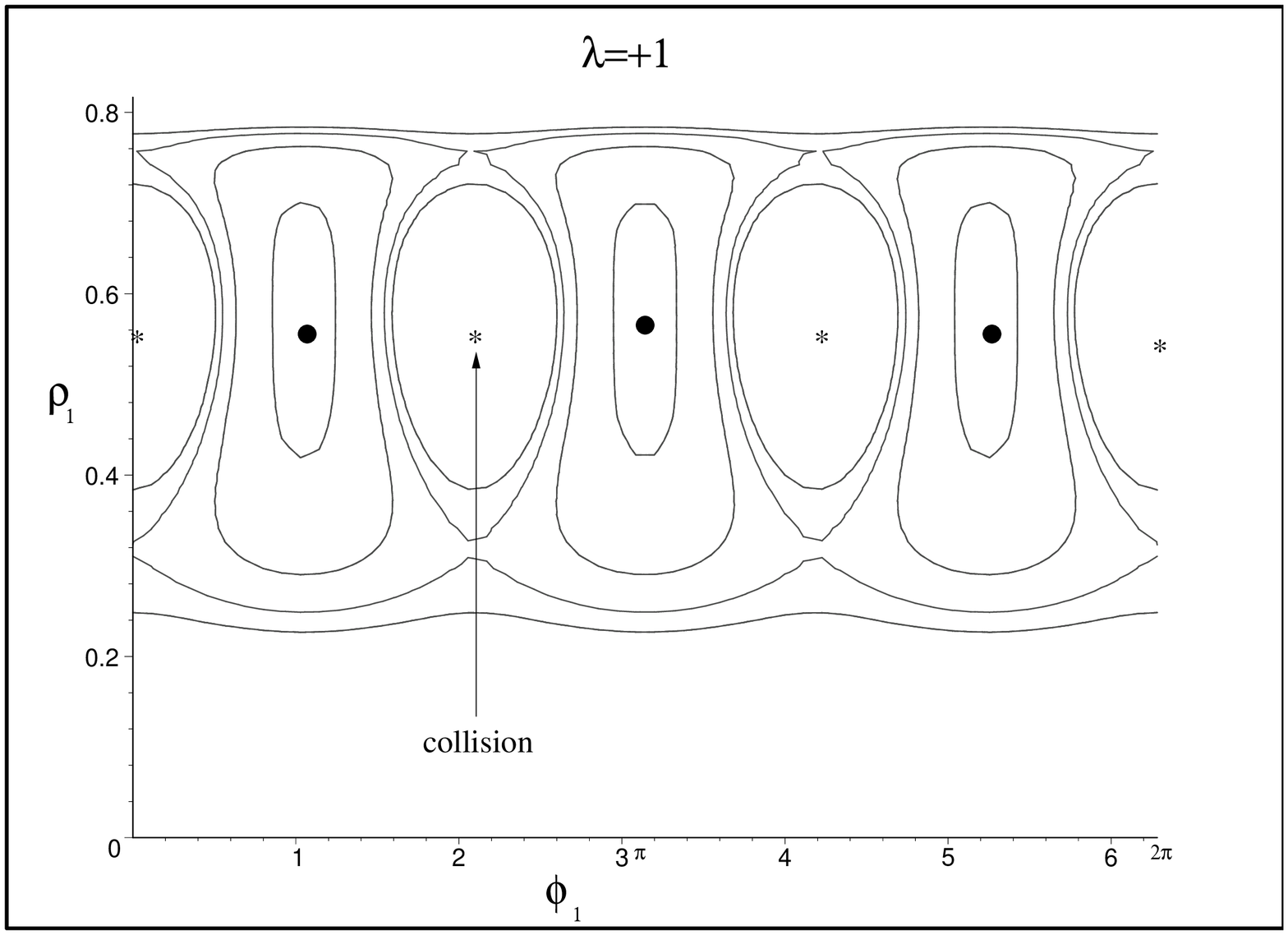}}
\caption{Configurations $C_N(2R,p)$. Phase portraits in $M_{\mu=1}$ for $N=3$ and
 different values of $\lambda$
(dots represent the trajectories for $\rho_{N+1}$, except for figure $(c)$ where the curves superpose).}
\label{planarfig}
\ec \efig

\bigskip

\noindent\textbf{\large Acknowledgements.}
The author warmly thanks James Montaldi, Tadashi Tokieda and  Anik Souli\`{e}re for helpful comments
and suggestions.

\bigskip

\noindent\textbf{\huge Appendix}

\bigskip

\noindent We describe in that section the possible types of point-orbits for each finite subgroups of
$SO(3)$ and $O(3)$.
The nomenclature of the two following tables is given in Section \ref{RPOident}.

An $n$-ring is a regular $n$-gon.  \pg Vertical\pd  refers to the axis of the rotation subgroup in
each case.  \pg Vertically aligned\pd means that the vortices of one ring
are directly above those of the other ring, and \pg vertically staggered\pd
means that the upper one is rotated by $\pi/n$ with respect to the
lower. Finally, $C'$ denotes the configuration dual to $C$.
\begin{table}[ht]   
$$\begin{array}{|c|c|c|c|l|}
\hline
& & & &  \\[-6pt]
\Gamma&  \kO&  I &  |\kO| & \mbox{Description} \\[4pt]
\hline\hline
\CC_n &  R &  1   & n  &  \mbox{$n$-ring} \\
    &  p &  \CC_n & 1  &   \mbox{pole} \\
\hline
\DD_n &  R &  1   & 2n &   \mbox{pair of $n$-rings on opposite
                                                        latitudes} \\
    &  r &  \CC_2 & n  &   \mbox{equatorial $n$-ring or dual} \\
    &  p &  \CC_n & 2  &    \mbox{pair of poles} \\
\hline
\Tet &  R &  1    & 12  &   \mbox{regular $\Tet$ orbit} \\
     &  e &  \CC_2  & 6   &    \mbox{mid-points of edges of tetrahedron} \\
     &  v &  \CC_3  & 4   &  \mbox{vertices of tetrahedron or dual} \\
\hline
\Oct &  R &  1    & 24  &   \mbox{regular $\Oct$ orbit} \\
     &  e &  \CC_2  & 12  &    \mbox{mid-points of edges of octahedron} \\
     &  f &  \CC_3  & 8   &    \mbox{mid-points of faces of octahedron} \\
     &  v &  \CC_4  & 6   &    \mbox{vertices of octahedron} \\
\hline
\Icos &  R &  1    & 60  &   \mbox{regular $\Icos$ orbit} \\
      &  e &  \CC_2  & 30  &    \mbox{mid-points of edges of icosahedron} \\
      &  f &  \CC_3  & 20  &    \mbox{mid-points of faces of icosahedron} \\
      &  v &  \CC_5  & 12  &    \mbox{vertices of icosahedron} \\
\hline
\end{array}$$
\caption{Classification of point-orbits of finite subgroups of $\SO(3)$.}
\label{so(3)table}
\end{table}

\begin{table}[ht]  
$$\begin{array}{|c|c|c|c|l|}
\hline
& & & &  \\[-6pt]
\Gamma&  \kO&  I &  |\kO| &  \mbox{Description} \\[4pt]
\hline\hline
\CC_{nv} &  R_{s} &  1      & 2n  &   \mbox{semi-regular $2n$-gon} \\
         &  R,\ R' &  \CC_h    & n   &  \mbox{regular $n$-ring or dual} \\
         &  p     &  \CC_{nv} & 1   &   \mbox{pole} \\
\hline
\CC_{nh} &  R     &  1      & 2n  & \mbox{pair of $n$-rings
                                               on opposite latitudes} \\
       &  R^e  &  \CC_h    & n   &  \mbox{equatorial $n$-ring} \\
       &  p      &  \CC_{n} & 2   &   \mbox{pair of poles} \\
\hline
\DD_{nh} &  R_{s} &  1      & 4n  &   \mbox{vertically aligned pair of
                                                  semi-regular $2n$-gons} \\
       &  R_{s}^e &  \CC_h    & 2n  &  \mbox{equatorial semi-regular
                                                                $2n$-gon} \\
       &  R,\ R'  &  \CC_h'  & 2n &  \mbox{vertically aligned pair
                                               of $n$-rings or duals} \\
       &  r,\ r'  &  \CC_{2v} & n   &   \mbox{equatorial
                                                       $n$-ring or dual} \\
       &  p      &  \CC_{nv} & 2   &  \mbox{pair of poles} \\
\hline
\DD_{nd} &  R_{s} & 1       & 4n  &   \mbox{vertically staggered pair of
                                                  semi-regular $2n$-gons}  \\
         &  R & \CC_h   & 2n  &  \mbox{vertically staggered pair of
                                                   $n$-rings} \\
         &  r & \CC_2   & 2n  &   \mbox{equatorial $2n$-ring}\\
         &  p      & \CC_{nv}& 2   &   \mbox{pair of poles} \\
\hline
\SS_{2n} &  R & 1       & 2n  &   \mbox{vertically staggered pair of
                                                   $n$-rings}   \\
         &  p      & \CC_n   & 2   &   \mbox{pair of poles} \\
\hline
\CC_h    &  R    &  1      & 2   &   \mbox{vertically aligned pair of points}  \\
         &  E      &  \CC_h  & 1   &  \mbox{equatorial point}\\
\hline
\CC_i    &  R    &  1      & 2   &   \mbox{pair of antipodal points}   \\
\hline
\Tet_d   &  R      &  1      & 24  &  \mbox{regular $\Tet_d$ orbit} \\
         &  E      & \CC_h   & 12  &  \mbox{generic orbit on edges of
                                                             tetrahedron} \\
         &  e      & \CC_{2v} & 6  &   \mbox{mid-points of edges of
                                                            tetrahedron} \\
         &  v,v^\prime      & \CC_{3v} & 4  &   \mbox{vertices of tetrahedron
                                                              or  dual}\\
\hline
\Tet_h   &  R      &  1      & 24  &   \mbox{regular $\Tet_h$ orbit}  \\
         &  E      &  \CC_h  & 12  &   \mbox{generic orbit on `equator'} \\
         &  e    & \CC_{2v}   & 6   &   \mbox{mid-points of edges of
                                                             tetrahedron} \\
         &  v      & \CC_3   & 8   &   \mbox{vertices of cube}\\
\hline
\Oct_h   &  R      &  1      & 48  &   \mbox{regular $\Oct_h$ orbit} \\
         &  E   &  \CC_h  & 24  &   \mbox{generic orbit on edges of
                                                             octahedron} \\
         &  E'  &  \CC_h'  & 24  &  \mbox{generic orbit on
                                             face bisectors of octahedron} \\
         & e    &  \CC_{2v}  & 12 &  \mbox{mid-points of edges of
                                                             octahedron}\\
         & f    &  \CC_{3v}  & 8  &   \mbox{mid-points of faces of
                                                             octahedron}\\
         & v    &  \CC_{4v}  & 6  &   \mbox{vertices of octahedron}\\
\hline
\Icos_h  &  R      &  1      & 120 &   \mbox{regular $\Icos_h$ orbit} \\
         &  E  &  \CC_h     & 60 &   \mbox{generic orbit on edges
                                               of icosahedron}\\
         & e    &  \CC_{2v}  & 30 &   \mbox{mid-points of edges of
                                                  icosahedron}\\
         & f    &  \CC_{3v}  & 20  &   \mbox{mid-points of faces of
                                                             icosahedron}\\
         & v    &  \CC_{5v}  & 12  & \mbox{vertices of icosahedron}\\
\hline
\end{array}$$
\caption{Classification of point-orbits of finite subgroups of $\OO(3)$.}
\label{o(3)table}
\end{table}


\bigskip

\begin{thebibliography}{99}
\bibitem[Ar82] {Ar82} H. Aref, Point vortex motion with a center of symmetry.
                      \emph{Phys. Fluids} {\bf 25} (1982), 2183-2187.
\bibitem[Ar83] {Ar83'} H. Aref, The equilibrium and stability of a row of point  vortices.
                       \emph{J. Fluid Mech.} {\bf 290} (1983), 167-181.
\bibitem[AV98] {Ar98} H. Aref and C. Vainchtein, Asymmetric equilibrium patterns of point vortices.
                       \emph{Nature} {\bf 392} (1998), 769-770.
\bibitem[ANSTV]{ANSTV03} H. Aref, P. Newton, M. Stremler, T. Tokieda, and DL. Vainchtein, Vortex crystals.
                       \emph{Preprint} (2002).
\bibitem[BC] {BC01} S. Boatto and H. Cabral, Non-linear stability of \ria of vortices on a
                       non-rotating sphere.  To appear in \emph{SIAM J. of Applied Math}.
\bibitem[CS99] {CS99} H. Cabral and D. Schmidt, Stability of \ria in the problem of $N+1$ vortices.
                       \emph{SIAM J. Math. Anal.} {\bf 31} (1999), 231-250.
\bibitem[H] {H} H. Helmhotz, On integrals of the hydrodynamical equations which express vortex motion,
                      \emph{Phil. Mag.} {\bf 33} (1867), 485-512.
\bibitem[KN98] {KN98} R. Kidambi and P. Newton, Motion of three point vortices on a sphere.
                       \emph{Physica D} {\bf 116} (1998), 143-175.
\bibitem[LMR01] {LMR00} C. Lim, J. Montaldi, M. Roberts, Relative equilibria of point vortices on the sphere.
                       \emph{Physica D} {\bf 148} (2001), 97-135.
\bibitem[LMR] {LMR} F. Laurent-Polz, J. Montaldi, M. Roberts, Stability of point vortices on the sphere.
                       Should be submitted to \emph{SIAM} soon.
\bibitem[LP02] {LP02} F. Laurent-Polz, Point vortices on the sphere: a case with opposite vorticities.
                       \emph{Nonlinearity} {\bf 15} (2002), 143--171.
\bibitem[LP] {LP} F. Laurent-Polz, Point vortices on a rotating sphere.
                       Preprint INLN (2002) (Available at arXiv: math.DS/0301360).
\bibitem[LPth] {LPTh} F. Laurent-Polz, Etude G\'{e}om\'{e}trique de la Dynamique de $N$ Tourbillons Ponctuels
                    sur une Sph\`{e}re. Ph.D. Thesis, University of Nice (2002).
\bibitem[LR96] {LR96} D. Lewis  and T. Ratiu, Rotating $n$-gon/$kn$-gon vortex configurations.
                       \emph{J. Nonlinear Sci.} {\bf 6} (1996), 385-414.
\bibitem[LT99] {LT99} L. Lerman  and T. Tokieda, On relative normal modes.
                         \emph{CR Acad. Sci. Paris}, Série I {\bf 328} (1999), 413-418.
\bibitem[Mo97] {Mo97b} J. Montaldi, Persistence d'orbites p\'{e}riodiques relatives dans les syst\`{e}mes
                        Hamiltonien sym\'{e}triques.
                       \emph{CR Acad. Sci. Paris}, S\'{e}rie I {\bf 324} (1997), 553-558.
\bibitem[MR94] {MR94} J. Marsden and T. Ratiu, Introduction to Mechanics and Symmetry.
                       TAM {\bf 17} Springer-Verlag (1994).
\bibitem[Or98] {Or98} J-P. Ortega, Symmetry, Reduction, and Stability in Hamiltonian Systems.
                        Ph.D. Thesis. University of California, Santa Cruz (1998).
\bibitem[OR] {OR} J-P. Ortega and T. Ratiu, The dynamics around stable Hamiltonian relative equilibria.
                        \emph{Preprint INLN} (2000).
\bibitem[P79] {P79} R. Palais, Principle of symmetric criticality.
                       \emph{Comm. Math. Phys.} {\bf 69} (1979), 19-30.
\bibitem[PM98] {PM98} S. Pekarsky and J. Marsden, Point vortices on a sphere: Stability of relative equilibria.
                       \emph{J. Math. Phys.} {\bf 39} (1998), 5894-5907.
\bibitem[ST] {ST} A. Souli\`{e}re and T. Tokieda, Periodic motions of vortices on surfaces with symmetry.
                         \emph{J.\ Fluid Mech.} {\bf 460} (2002) 83--92.
\bibitem[To01] {To} T. Tokieda, Tourbillons dansants.
                         \emph{CR Acad. Sci. Paris}, S\'{e}rie I {\bf 333} (2001), 943-946.
\bibitem[WR] {WR} C. Wulff and M. Roberts, Hamiltonian systems near relative periodic orbits.
                       To appear in SIAM.
\end{thebibliography}
\end{document}